\documentclass[12pt]{article}
\usepackage{amsmath}
\usepackage{graphicx,psfrag,epsf}
\usepackage{enumerate}
\usepackage[comma,authoryear]{natbib}%
\usepackage{url} 

\usepackage{amsthm,amsmath}
\usepackage{graphicx}
\usepackage{amssymb}
\usepackage[utf8]{inputenc}
\usepackage[mathscr]{eucal}
\usepackage{fancyhdr}
\usepackage{fancybox}
\usepackage{color}
\usepackage{longtable}
\usepackage{xcolor}
\usepackage{paralist}
\usepackage{undertilde}
\usepackage{bbm}
\usepackage{paralist}

\usepackage{stackengine}
\stackMath
\newcommand\tenq[2][1]{%
\def\useanchorwidth{T}%
\ifnum#1>1%
\stackunder[0pt]{\tenq[\numexpr#1-1\relax]{#2}}{\scriptscriptstyle\thicksim}%
\else%
\stackunder[1pt]{#2}{\scriptscriptstyle\thicksim}%
\fi%
}

\RequirePackage[colorlinks,citecolor=blue,urlcolor=blue]{hyperref}

\newcommand{\pmss}{{\pmss}}
\newcommand{\pr}{^{\prime}}
\newcommand{\n}{^{(n)}}

\newcommand{\npr}{^{(n)\prime}}

\newcommand{\E}{\mathrm{E}}

\newcommand{\R}{\mathbbm{R}}

\newcommand{\pms}{{\scriptscriptstyle \pm}}

\newcommand{\Fs}{\mathbf{F}_{\pms}^{(n)}}

\newcommand{\tZ}{\mathbf{Z}}

\newcommand{\barc}{\bar{c}}
\newcommand{\tu}{\mathbf{u}}

\newcommand{\tT}{\mathbf{T}}

\newcommand{\td}{\mathbf{d}}

\newcommand{\tmu}{\boldsymbol{\mu}}
\newcommand{\tSigma}{\boldsymbol{\Sigma}}

\renewcommand{\phi}{\varphi}

\newtheorem {Proposition}{Proposition} [section]

\newtheorem {Corollary}[Proposition] {Corollary}

\newtheorem {Assumption}{Assumption}[section]

\newtheorem*{Proof*}{Proof}

\usepackage[splitrule]{footmisc}
\usepackage{xr}

\begin{document}

\def\spacingset#1{\renewcommand{\baselinestretch}%
{#1}\small\normalsize} \spacingset{1}


{
  \title{\bf Efficient Fully Distribution-Free \\ Center-Outward Rank Tests \\  for Multiple-Output Regression\\ and MANOVA}
  \author{Marc Hallin
   \hspace{.2cm}\\
    ECARES and  D\' epartement de Math\' ematique\\ Universit\' e libre de Bruxelles, Brussels, Belgium\medskip\\
    Daniel Hlubinka and  \v S\' arka Hudecov\' a 
    \hspace{.2cm}\ \\
    Faculty of Mathematics and Physics\\ Charles University, Prague,  Czech Republic\vspace{-1mm}
   }
     \date{}
  \maketitle



\begin{abstract}
Extending rank-based inference to a multivariate setting such as multiple-output regression or MANOVA  with unspecified $d$-dimen\-sional error density has remained an open problem for more than half a century. None of the many solutions proposed so far is enjoying the combination of distribution-freeness and efficiency that makes rank-based inference a successful tool in the univariate setting. A concept of {\it center-outward} multivariate ranks and signs based on measure transportation ideas has been introduced recently. Center-outward ranks and signs are not only distribution-free but achieve in dimension $d>1$ the (essential) maximal ancillarity property of traditional univariate ranks. In the present case, we show that  fully distribution-free testing procedures based on center-outward ranks can achieve parametric efficiency. We establish  the H\' ajek representation and asymptotic normality results required in  the construction of such tests in multiple-output regression and MANOVA models. Simulations and an empirical study demonstrate the excellent performance of the proposed procedures. 
\end{abstract}
\noindent%
\vspace{-5mm}

\spacingset{1.45} 

\section{Introduction} Linear models---regression (single- and multiple-output), Analysis of Variance (ANOVA and MANOVA)---are probably the most popular and most useful of all statistical models; they are found in the table of contents of all statistical textbooks and  statistical softwares, and  are part of daily statistical practice in all domains of application. The pseudo-Gaussian approach---Gaussian quasi maximum likelihood estimation and pseudo-Gaussian~$F$ tests---is largely dominant in that context,  
on the ground that pseudo-Gaussian methods   remain asymptotically valid under a broad class of non-Gaussian densities satisfying mild moment conditions. One should beware of excessive confidence in such asymptotics, 
though. \vspace{-2mm}

\subsection{Pseudo-Gaussian tests}\label{pseudoSec} Let us concentrate on hypothesis testing. The problem with pseudo-Gaussian tests under unspecified noise density is twofold: 
\begin{compactenum}
\item[(a)]although pseudo-Gaussian tests are asymptotically valid under a broad range of non-Gaussian densities, that asymptotic validity is far from uniform: actually,  in a semiparametric  model with parameter $\theta$ where the underlying noise has unspecified density~$f$ in some broad class $\cal F$ of densities, 
 a sequence $\phi^{(n)}$ of  tests of  the null hypothesis $\theta = \theta_0$ has asymptotic level $\alpha$ iff $\lim_{n\to\infty}\sup_{f\in{\cal F}}{\rm E}_{\theta_0,f}[\phi^{(n)}]\leq\alpha$, whereas pseudo-Gaussian tests~$\phi^{(n)}_{\cal G}$ only satisfy the pointwise condition $\lim_{n\to\infty}{\rm E}_{\theta_0,f}[\phi^{(n)}_{\cal G}]\leq\alpha$ for all $f\in{\cal F}$; 
\item[(b)]still for fixed $n$, the performance of pseudo-Gaussian tests may rapidly deteriorate away from the Gaussian.%
\end{compactenum}
Appendix~A.1 illustrates these pitfalls in the case of Hotelling's bivariate two-sample test.

\subsection{Rank-based tests} A natural way to restore uniform asymptotics, thereby solving the validity problem in~(a) consists in  resorting to distribution-free tests, and this is how rank tests  enter the picture.  
 Rank-based testing methods have been quite
successful in testing problems for  
single-ouput regression and linear models such as ANOVA (see the
classical monographs by \cite{HajekSidak},
\cite{RWolfe79}  
or \cite{PuriSen1985}) and univariate linear time series
(\cite{Hallinetal1985}, \cite{KoulSaleh1993},
\cite{HallinPuri1994}).  
Being distribution-free, rank tests remain valid over the full class of
absolutely continuous distributions.  In linear models (this includes
testing for single-output regression slopes, testing for treatment
effects in analysis of variance, testing against location shifts in
two-sample problems) and ARMA time series, they do reach parametric or
semiparametric efficiency bounds at given reference densities, thus
reconciling the conflicting objectives of robustness and efficiency. 
The celebrated Chernoff-Savage result (\cite{chersav58} and, for time-series, \cite{hal94}) moreover indicates that, far from losing power with respect to their pseudo-Gaussian counterparts, rank tests strictly dominate them under any non-Gaussian  density  $f$, making the latter non-admissible.

Extending these attractive features to a multivariate (multiple-output) context, of course, is highly desirable and the problem of defining multivariate concepts of ranks has been
a long-standing open problem, for which many solutions have been
proposed in the literature.  \cite{PuriSen1971} for a variety of
problems in multivariate analysis (including multiple-output
regression and MANOVA) and \cite{HallinIngPuri1989} for VARMA time
series models construct tests based on {it componentwise ranks} which,
however, fail to be distribution-free.  Building upon an ingenious
multivariate extension of the L$_1$ definition of quantiles,
\cite{Oja1999,Oja2010}  
defines the so-called {\it spatial ranks}; the resulting tests are
neither distribution-free nor efficient. Tests based on the ranks of
various concepts of statistical depth also have been proposed
(\cite{Liu92}, \cite{LiuSingh93},   \cite{ZuoHe06}).
While distribution-free, these ranks   are
failing to exploit any directional information, and hence typically do
not allow for any type of asymptotic efficiency. 
     As for the tests based on the {\it Mahalanobis ranks
  and signs} proposed by \cite{HallinPain2002b, HallinPain2002a,
  HallinPain2004, HallinPain2005},  
they do achieve, within the class of linear models and linear time
series with elliptical densities, parametric or semiparametric
efficiency at correctly specified elliptical reference densities;
their distribution-freeness, hence their validity, unfortunately, is
limited to the class of elliptical distributions. 

Inspired by measure transportation ideas, a new concept of ranks and
signs for multivariate observations has been introduced recently under
the name of {\it Monge-Kantorovich ranks and signs} in
\cite{Chernozhukov},  
under the name of {\it center-outward ranks and signs} in
\cite{Hallin2017} and \cite{Hallinetal2020}, along with the related population 
concepts of {\it center-outward distribution and quantile
  functions}. Unlike earlier
concepts, these ranks and signs extend to dimension~$d>1$ the
{\it essential maximal ancillarity} property (see Section~2.4~and Appendices D1 and D.2 of \cite{Hallinetal2020})    of univariate ranks; 
the corresponding empirical center-outward distribution functions,
moreover,  
satisfy a Glivenko-Cantelli result.

Center-outward ranks and signs  have been successfully applied  (\cite{Boeckel2018MultivariateBC}, \cite{BodhiDeb, BodhiGhos, sdhh21, Shietal, shi2020rateoptimality}) in the construction of distribution-free tests  of   independence between random vectors and multivariate goodness-of-fit; applications to the study of tail behavior and extremes can be found  in \cite{DeValk2018};  \cite{Beirlantetal2020} are using the related center-outward empirical quantiles in the analysis of multivariate risk;  \cite{HLL19, HLL20} are proposing  center-outward tests and R-estimators for VAR and VARMA time series models with unspecified innovation densities.  We refer to \cite{H21} for a review. The present paper goes one step further in the direction of a toolkit of distribution-free tests for multiple-output multivariate analysis 
by deriving a H\' ajek-type asymptotic representation result for 
linear center-outward rank statistics. Asymptotic normality   follows as a corollary, from which
center-outward rank  tests are constructed for multiple-output
regression models (including, as special cases, MANOVA and two-sample
location models). Those tests are fully distribution-free, hence valid,  over the entire
family of absolutely continuous distributions;  for adequate choice of the scores, parametric   
efficiency is attained at chosen  densities. 
 Since this paper was written  \citep{HHH2020}, some further results (among them,  partial Chernoff-Savage and Hodges-Lehmann properties) on the particular case of the two-sample location problem have been obtained  by \cite{deb2021}; see \cite{hm21} for some numerical comparisons with the tests presented here.

\subsection{A motivating example}\label{motivating}

The  importance of  center-outward rank tests in daily statistical practice is illustrated with the following real-life motivating example.   The Wisconsin Diagnostic Breast Cancer  data (WDBC; dataset  available at Machine Learning Repository \cite{Dua2019}), first analyzed in  \cite{street} in a classification context,  contains records on $n=569$ patients from two groups---benign or malignant tumor diagnosis.  
For each patient, several features were recorded from the digitized image of a fine needle aspirate of the breast mass, resulting in $d=30$ variables, labeled V1--V30. The two groups of patients are well separated: the two-sample Hotelling test  in dimension~$d=~\!30$   very significantly  rejects the null hypothesis of equal locations  (the R program delivers a $p$-value 0.000, meaning that the actual  $p$-value is less than $10^{-22}$!).{So does the Wilcoxon center-outward rank test. }

{In real life, however, diagnoses requiring 30 clinical measurements are highly impractical and costly. Reducing that number from 30 to 3 or 4 without losing diagnostic efficiency is an important issue. 
 Unfortunately,   when restricted to subsets of three or four variables,  the Hotelling test  typically  is   inconclusive}. Consider, for instance,  the subset consisting of  V12  (mean of  fractal dimension), V14 (standard error  of texture),  
 V21 (standard error  of symmetry),  and V22 (standard error of fractal dimension).   
   Figure \ref{fig-wdbc} shows bivariate scatterplots and histograms for  these four variables, revealing 
    skewness in all  univariate marginals and deviations from  elliptical symmetry in bivariate marginals.

 \begin{figure}[ht!]
\centering
\includegraphics[width=0.8\textwidth]{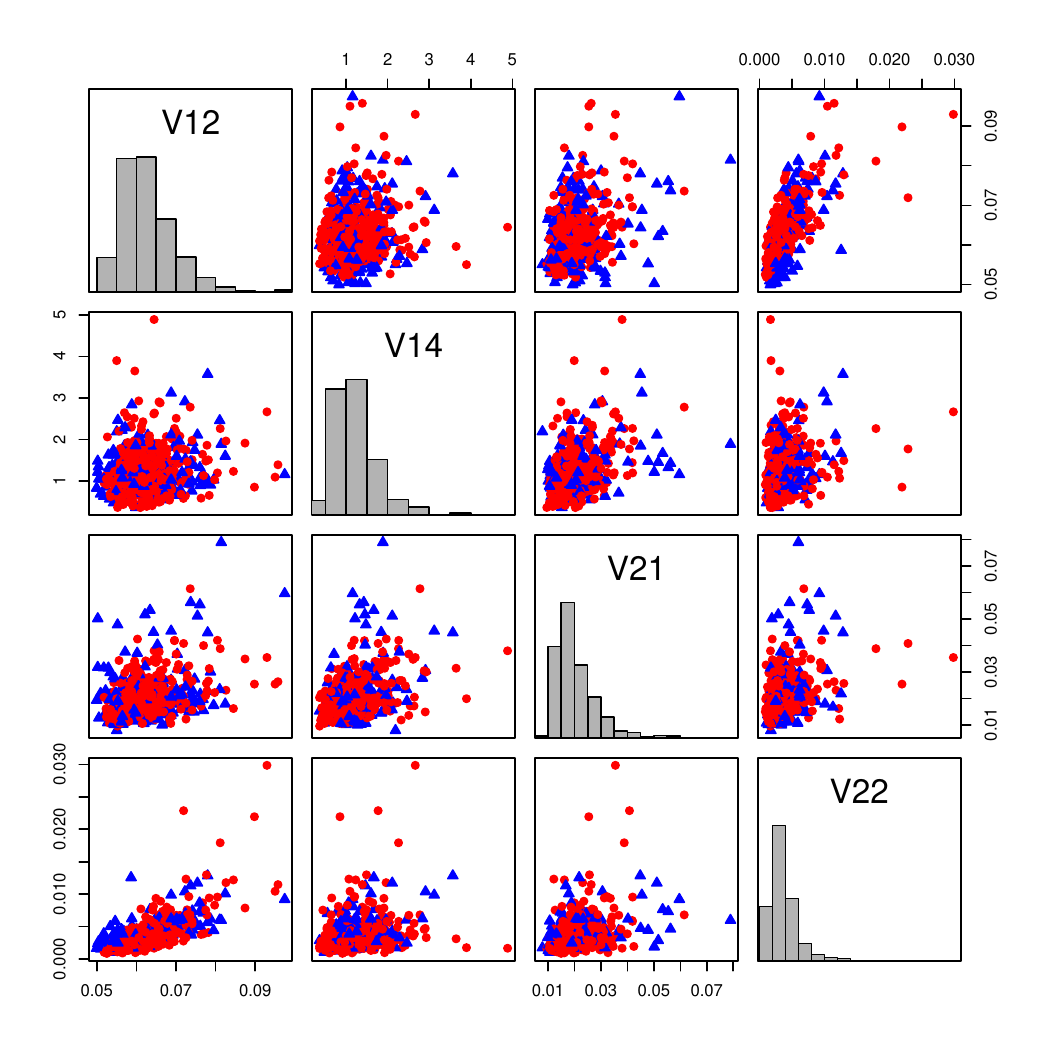}\vspace{-7mm}
 \caption{\small   
 Wisconsin Diagnostic Breast Cancer (WDBC) data: bivariate scatterplots\vspace{-2mm} and univariate histograms for mean fractal dimension (V12), standard error of texture\vspace{-2mm} (V14), 
  standard error of symmetry (V21),  and standard error of fractal dimension (V22) in 212 malignant\vspace{-2mm} patients (triangles) and 357 benign patients (circles).\vspace{-2mm}}\label{fig-wdbc}
\end{figure}

The Hotelling  
 and Wilcoxon  center-outward rank tests (see Section~\ref{sec-2slt} for a precise description) have been performed for the corresponding  four-dimensional  dataset and all its three-dimensional margin\-als\footnote{The  center-outward ranks were computed for a $569$-point random grid\vspace{-2mm} 
with  $n_R=20$, $n_S=28$, $n_0=9$; the~$n_S=28$-points  over the sphere were generated  (seed \texttt{1111} in \cite{R})  
 as   in Section \ref{sec_d6}.}; $p$-values are shown in Table~\ref{tab3}. 
 With $p$-value 0.0090, the Wilcoxon test in dimension~4 is significant  at 5\% and 1\% levels, while Hotelling \linebreak (with $p$-value~0.0595) is not. Turning to dimension 3, Wilcoxon is always significant at~5\% level (at 1\% level in all cases but one), while Hotelling  never rejects on 1\%. The most spectacular case  is that of the subset $\{V12, V14, V21\}$ where Wilcoxon and Hotelling yield $p$-values  0.0327 and  0.9899, respectively.  
  \begin{table}[t!]
\centering
\begin{tabular}{cccccc}
  \hline  \hline 
 Variables & (12,14,21, 22)& (12,14,21) &(12,14,22) &  (12,21,22) & (14,21,22) \\ 
   \hline
   Hotelling & 0.0595 & 0.9899 & 0.0299 &  0.0346 &  0.2136\\
  c-o Wilcoxon & 0.0090&  0.0327& 0.0007&  0.0000 & 0.0018  \\ 
   \hline   \hline
\end{tabular}
\caption{\small Wisconsin Diagnostic Breast Cancer (WDBC) data:\vspace{-2mm} $p$-values of the two-sample location  Hotelling and center-outward Wilcoxon\vspace{-2mm} rank tests  
for the 4-dimensional marginal WDBC data corresponding to the set of variables $\{$V12,V14,VV21,V22$\}$ and its~3-dimensional subsets.  
}\label{tab3}
\end{table}
Such discrepancies most likely originate in the skewness, non-ellipticity and/or the heavy tails of the observations; their impact in terms of diagnostic power   may have crucial consequences.

\subsection{Outline of the paper} The paper is organized as follows. Section 2 briefly describes the main tools to be used: center-outward distribution and quantile functions (Section~2.1) and their empirical counterparts, the center-outward ranks and signs  (Section 2.2). The main properties of these concepts are summarized in Section 2.3 (Proposition 2.1); their invariance/equivariance properties are established in Proposition 2.2. Section 3 is entirely devoted to the key theoretical results of this paper, which extend and generalize the classical approach by H\' ajek and \v Sid\' ak~(1967): a H\' ajek-type asymptotic representation for multivariate center-outward linear rank statistics  and the resulting asymptotic normality result. Section 4.1 describes the multiple-output regression model to be considered throughout, which contains, as particular cases, the two-sample location and MANOVA models, of obvious practical importance. Local asymptotic normality is established in Section~4.2 for this model under general error densities (Proposition~4.1) and, for the purpose of future comparisons, for the particular case of elliptical distributions (Proposition~4.2). The center-outward rank tests we are proposing are described in Section~5.2, along   with (Corollary~5.2) their local   asymptotic optimality properties. Due to their importance in applications, the particular cases of the hypotheses of equal locations in the two-sample problem and  no treatment effect in MANOVA are considered in Section~5.3. Sections~6.1 and~6.2 propose some simple choices of score functions, extending the classical median-test-score (based on center-outward signs only), Wilcoxon, and van der Waerden (normal-score) tests. Section~6.3 discusses affine invariance issues. Section~7 is devoted to a Monte Carlo exploration, in dimension $d=2$, of the finite-sample performance of our rank tests which appear to outperform their competitors in non-elliptical situations  while performing equally well under ellipticity. Section~7.3 presents an archaeological  MANOVA application  in dimension~$d=4$;  
 while traditional MANOVA methods cannot reject the hypothesis of no treatment effect, our fully distribution-free center-outward  rank-based test rejects it quite significantly, which might lead to revising some of  the conclusions   \citep{data-source} on Middle-East
 economic exchanges between Egypt and Syro-Palestine in the {B}yzantine-{I}slamic transition period.  All proofs are concentrated in an online appendix where we also provide simulations in~dimension~$d=6$. 
 \setcounter{equation}{0}

\section{Center-outward    ranks and signs in~$\mathbb{R}^d$}\label{defsec}

\subsection{Center-outward distribution 
  functions}\label{FQsec}
   Throughout, denote  by~$\tZ\n$  a
triangular array $(\tZ_1\n,\dots, \tZ_n\n)$, $n\in\mathbb{N}$ of i.i.d.\ $d$-dimensional random vectors with distribution $\mathrm{P}$ in the family ${\mathcal P}_d$ of 
absolutely
continuous distributions on $\R^d$.  The notation $\overline{\text{spt}}({\rm P})$ is used for the support of ${\rm P}$,  ${\text{spt}}({\rm P})$ for its interior. The open (resp. closed) unit ball and the
unit hypersphere in ${\mathbb R}^d$ 
 are denoted  by~${\mathbb{S}_d}$ (resp. $\overline{{\mathbb{S}}}_d$) 
and~${\mathcal{S}_{d-1}}$, respectively; ${\rm U}_d$ stands for the spherical\footnote{Namely, the spherical\vspace{-2mm}
distribution with uniform (over $[0,1]$) radial density---equivalently,
the product of a uniform over the distances to the origin and a
uniform over the unit sphere ${\cal S}_{d-1}$.\vspace{-2mm} For $d=1$, it coincides with the Lebesgue uniform; for $d\geq 2$, it has unbounded density at the origin.} uniform distribution over  ${\mathbb{S}_d}$, $\mu_d$   for the Lebesgue measure over~$\mathbb{R}^d$; ${\bf I}_d$ is the $d\times d$ unit matrix, ${\bf 1}_A$ the indicator of the Borel set $A$.  

The definition of the center-outward distribution function of $\rm P$ is particularly simple for~${\rm  P}$ in the so-called class ${\cal P}^+_d$ of distributions {\it with nonvanishing densities}---namely, the class of all distributions with density $f:={\rm d P}/{\rm d}\mu_d$ such that, for all $D\in\mathbb{R}^+$, there exist constants~$\lambda^-_{D;\mathrm{P}}$ and  $\lambda^+_{D;\mathrm{P}}$ satisfying 
 $ 0<\lambda^-_{D;\mathrm{P}}\leq f({\bf z})
\leq \lambda^+_{D;\mathrm{P}}<\infty$   for all $\bf z$ with~$\Vert{\bf z}\Vert \leq D$ (so that~spt$({\rm P})=\mathbb{R}^d$ and $\rm P$-a.s.\ is equivalent to $\mu_d$-a.e.). The main result in \cite{McCann} then implies the existence of an a.e.~unique convex lower semi-continuous function~$\phi:~\!\mathbb{R}^d\to~\!\mathbb{R}$ with  gradient $\nabla\phi$ such that~$\nabla\phi\#{\rm P}={\rm U}_d$---we borrow from measure transportation   the convenient  notation  $T\#\mathrm{P}$ ($T\!:\mathbb{R}^d\to~\!\mathbb{R}^d$   {\it pushes~$\mathrm{P}$ forward to~$T\#\mathrm{P}$}) for the distribution  under~${\bf Z}\sim\mathrm{P}$ of~$T({\bf Z})$.  Call~${\bf F}_\pms:=\nabla\phi$ the {\it center-outward distribution function} of $\rm P$. It follows from \cite{Figalli} that ${\bf F}_\pms$ defines a homeo\-morphism  between the punctured unit ball~${\mathbb S}_d\!\setminus\!\{{\bf 0}\}$ and its image~$\mathbb{R}^d\!\setminus\!{\bf F}_\pms^{-1}({\bf 0})$:
   call~${\bf Q}_\pms: {\bf u}\mapsto {\bf Q}_\pms({\bf u}):= {\bf F}_\pms^{-1}({\bf u})$, ${\bf u}\neq{\bf 0}$ the  {\it center-outward quantile function}.   \cite{Figalli} also shows that, defining   ${\bf Q}_\pms({\bf 0}):= {\bf F}_\pms^{-1}({\bf 0})$ yields a convex and compact subset  with Lebesgue measure zero in $\mathbb{R}^d$, the  {\it center-outward median set  of} $\rm P$. 
 
 All the intuition and all the properties of center-outward distribution and quantile functions hold for 
 ${\rm P}\in{\cal P}^+_d$; this special case is the one considered in~\cite{Hallin2017}.  
  A more general case  
  is addressed in~\cite{BarrioGonzHallin} and \cite{Hallinetal2020} where we refer to for details,  
  but requires   more technical 
  definitions, which we are skipping here. 
  Note, however, that while some  statements below  only hold under   ${\rm P}\in{\cal P}_d^+$,  
  many others (including validity), due to distribution-freeness, can be made  under the 
  very general condition ${\rm P}\in{\cal P}_d$.

\subsection{Center-outward ranks and signs}\label{RSsec}

Except for a few particular cases such as spherical distributions,  the above definitions  are not meant for an analytical derivation of ${\bf F}_{\pms}$ and~${\bf Q}_{\pms}$ which typically involves Monge-Amp\`ere  
 equations; in particular, no closed forms of  ${\bf F}_{\pms}$ and ${\bf Q}_{\pms}$ are known for   non-spherical elliptical distributions.  Estimation is possible, though, via their empirical counterparts~${\bf F}_\pms\n$ and~${\bf Q}_\pms\n$, based on center-outward ranks and signs, which we now describe. 

Associated with the $n$-tuple $\tZ_1\n,\dots, \tZ_n\n$,  the {\it
  empirical center-outward distribution function} ${\bf F}_\pms\n$ is mapping
$\tZ_1\n,\dots, \tZ_n\n$ to a ``regular'' grid $\mathfrak{G}_n$ of the unit 
ball~${\mathbb S}_d$. That grid $\mathfrak{G}_n$ is obtained as follows:\smallskip 

\begin{compactenum}
\item[{\it (a)}] first factorize $n$ into $n=n_Rn_S + n_0$, with
  $0\leq n_0<\min(n_R, n_S)$;
\item[{\it (b)}] next consider a ``regular array"  
  $\mathfrak{S}_{n_S}:=\{{\bf s}^{n_S}_1,\ldots,{\bf s}^{n_S}_{n_S}\}$ of $n_S$ points on the sphere 
  ${\cal S}_{d-1}$ (see the comment below); 
\item[{\it (c)}] finally, the grid consists in the collection $\mathfrak{G}_n$ of  the 
  $n_Rn_S  $ points $\mathfrak{g}$ of the
  form\vspace{-1mm}
  \[
  \big(r/\big(n_R +1\big)\big){\bf s}^{n_S}_s,\quad 
  r=1,\ldots,n_R,\ \  s=1,\ldots,n_S,
\vspace{-1mm}  \]
   along with ($n_0$ copies of)  the
  origin in case $n_0\neq 0$: a total number $n-(n_0 -1)$ or $n$ of distinct points, thus, according as $n_0>0$ or $n_0=0$.
\end{compactenum}
By ``regular'' we mean ``as uniform as
  possible'', in the sense, for example, of the {\it
    low-discrepancy sequences} of the type considered in numerical
  integration and Monte-Carlo methods (see, e.g., \cite{Niederreiter}, \cite{Judd},  
   or \cite{Santner}).
    The only mathematical requirement needed for Proposition~\ref{H2018} below is the weak convergence, as~$n\to\infty$, of the uniform discrete distribution over~$\mathfrak{S}_{n}$  to the uniform distribution over~${\mathbb S}_{d}$;  all sequences~$\mathfrak{S}_{n}$ satisfying that requirement  yield the same asymptotic results. A uniform i.i.d.~sample of points over~${\mathbb S}_{d}$, for example, satisfies the requirement but fails to produce mutually independent  ranks and signs; moreover,  one easily can  construct arrays that are ``more regular"  than an i.i.d.~one.  For instance, one could  see that $n_S$ or $n_S-1$  
   of the points~${\bf s}^{n_S}_s$ in~$\mathfrak{S}_n$  are such that $-\, {\bf s}^{n_S}_s$ also belongs to~$\mathfrak{S}_{n_S}$, so that~$\Vert \sum_{s=1}^{n_S}{\bf s}^{n_S}_s\Vert $ is~0 or~1 according as $n_S$ is even or odd.  One also could consider factori\-zations of the\linebreak form $n=n_Rn_S + n_0$ with $n_S$ even and~$0\leq n_0<\min(2n_R, n_S)$, then require~$\mathfrak{S}_n$ to be symmetric with respect to the origin, automatically yielding~$\sum_{s=1}^{n_S}{\bf s}^{n_S}_s={\bf 0}$.

The empirical counterpart ${\bf F}_\pms\n$ of ${\bf F}_\pms$ 
 is defined as the (bijective, once the origin is given multiplicity $n_0$)  
mapping from~$\tZ_1\n,\dots, \tZ_n\n$ to the grid~$\mathfrak{G}_n$ that mini\-mizes the sum of squared Euclidean distances~$\sum_{i=1}^n\big\Vert {\bf F}_\pms\n (\tZ_i\n) - \tZ_i\n \big\Vert ^2$. That
mapping is unique with probability one; in practice, it is obtained
via a simple optimal assignment (pairing) algorithm (a linear program; see Section~4 of \cite{Hallin2017}  
for details).

Call {\it center-outward rank} of $\tZ_i\n$ the integer (in~$\{1,\ldots , n_R\}$ or~$\{0, \ldots , n_R\}$ according as~$n_0=0$ or not)
$R\n_{i;{{\pms}s}}:=(n_R +1)\big\Vert {\bf F}_\pms\n (\tZ_i\n)\big\Vert$ 
and {\it center-outward sign}  of~$\tZ_i\n$ the  unit vector
${\bf S}\n_{i;{\pms}}:={\bf F}_\pms\n (\tZ_i\n)/\big\Vert {\bf F}_\pms\n (\tZ_i\n)\big\Vert$ for ${\bf F}_\pms\n (\tZ_i\n)\neq{\bf 0}$; 
for ${\bf F}_\pms\n (\tZ_i\n) ={\bf 0}$, put~${\bf S}\n_{i;{{\pms}s}}={\bf 0}$. 

Some desirable finite-sample properties, such as strict independence between the ranks and the signs,  only hold for $n_0=0$ or 1, due to the fact that the mapping from the sample to the grid is no longer injective for $n_0\geq 2$. This, which  has no asymptotic consequences (since  the number $n_0$ of tied values involved is $o(n)$ as $n\to\infty$),    is easily taken care of by the following tie-breaking device: 
\begin{compactenum}
\item[{\it (i)}] randomly select  $n_0$ directions ${\bf s}^0_1,\ldots,{\bf s}^0_{n_0}$   in~$\mathfrak{S}_{n_S}$, then  
\item[{\it (ii)}] replace the $n_0$ copies  of the origin with the new gridpoints  $\frac{1}{2(n_R+1)}{\bf s}^0_1,\ldots,\frac{1}{2(n_R+1)}{\bf s}^0_{n_0}$. 
\end{compactenum}
The resulting grid (for simplicity, the same notation ${\mathfrak{G}}_n$ is used) no longer has multiple points, and the optimal pairing between the sample and the grid is bijective; the $n_0$ smallest ranks, however,  take the non-integer value~$1/2$. Again, this tie-breaking device has no influence on asymptotic results.\vspace{-0mm}

\subsection{Main properties}\label{Propsec} This section summarizes 
 the main properties of the concepts defined in Sections~\ref{FQsec} and~\ref{RSsec}; further properties and a proof for Proposition~\ref{H2018} can be found in \cite{Hallinetal2020}.

\begin{Proposition}\label{H2018} Let ${\bf F} _{\pms}$ denote the center-outward distribution function of~${\rm P}\in{\cal P}_d$.  Then,
\begin{compactenum}
\item[{\it (i)}]${\bf F} _{\pms}$ is a probability integral transformation of $\mathbb{R}^d$: namely, ${\bf Z}\sim {\rm P}$ iff~${\bf F} _{\pms}({\bf Z})\sim {\rm U}_d$; by construction, $\Vert{\bf F} _{\pms}({\bf Z})\Vert$ is uniform over the interval $[0, 1]$, ${\bf F} _{\pms}({\bf Z})/\Vert{\bf F} _{\pms}({\bf Z})\Vert$ uniform over the sphere ${\cal S}_{d-1}$, and they are mutually independent. 
\end{compactenum}
 Let ${\bf Z}\n_i,\ldots ,{\bf Z}\n_i$ be i.i.d.\ with distribution ${\rm P}\in{\mathcal P}_d$ and center-outward distribution function~${\bf F} _{\pms}$. Then,
 \begin{compactenum}
\item[{\it (ii)}] $\big({\bf F}\n_{\pms}({\bf Z}\n_{1}),\ldots , {\bf F}\n_{\pms}({\bf Z}\n_{n}) \big)$ is uniformly distributed over the $n!/n_0!$ permutations with repetitions  of  the  gridpoints in $\mathfrak{G}_n$ with the origin  counted as $n_0$ indistinguishable points (resp. the $n!$ permutations of $\mathfrak{G}_n$ if either $n_0\leq 1$ or the tie-breaking device described in Section~\ref{RSsec} is adopted);
           \item[{\it (iii)}]  if either $n_0=0$ or  the tie-breaking device described in Section~\ref{RSsec} is adopted, the $n$-tuple  of center-outward ranks $\big(R\n_{1;\pms }, \ldots , R\n_{n;\pms }\big)$ and the $n$-tuple of~center-out\-ward  signs~$\big({\bf S}\n_{1;\pms }, \ldots , {\bf S}\n_{n;\pms }\big)$ are mutually independent;
           \item[{\it (iv)}]  if either $n_0\leq 1$ or  the tie-breaking device described in Section~\ref{RSsec} is adopted, \linebreak the~$n$-tuple~$\big({\bf F}\n_{\pms}({\bf Z}\n_{1}),\ldots , {\bf F}\n_{\pms}({\bf Z}\n_{n}) \big)$ is {\em essentially maximal ancillary}.\footnote{See Section~2.4 and Appendices D1 and D.2\vspace{-2mm} of \cite{Hallinetal2020} for a precise definition of this crucial  property (which entails distribution-freeness) and a proof.}
           \vspace{-1mm}
\end{compactenum}
Assuming, moreover, that ${\rm P}\in{\mathcal P}_d^{+}$,
\begin{compactenum}
\item[{\it (v)}] (Glivenko-Cantelli) 
 $\displaystyle{\max_{1\leq i\leq n}}\Big\Vert   {\bf F}\n _{\pms}({\bf Z}\n_i) - {\bf F} _{\pms}({\bf Z}\n_i) \Big\Vert  \rightarrow 0$ \textrm{\em a.s.}  as~$n\to\infty$.
 \end{compactenum}
\end{Proposition}

Center-outward distribution functions, ranks, and signs also inherit, from the invariance features  of  Euclidean distances, elementary but quite remarkable invariance and equivariance properties under orthogonal transformations. Denote by ${\bf F}^{{\bf Z}}_{\pms}$ the center-outward distribution function of $\bf Z$ and by~${\bf F}^{{\bf Z};(n)}_{\pms}$ the empirical distribution function of a sample~${\bf Z}_1,\ldots,{\bf Z}_n$ associated with a grid $\mathfrak{G}_n$. 

\begin{Proposition}\label{invF} Let $\boldsymbol{\mu}\in\mathbb{R}^d$ and denote by ${\bf O}$ a $d\times d$ orthogonal matrix. Then,
\begin{compactenum}
   \item[{\it (i)}]  ${\bf F}^{\boldsymbol{\mu}+ {\bf O}{\bf Z}}_{\pms} (\boldsymbol{\mu} + {\bf O}{\bf z})= {\bf O}{\bf F}^{\bf Z}_{\pms}({\bf z})$, ${\bf z}\in\mathbb{R}^d$;
      \item[{\it (ii)}]  denoting by ${\bf F}^{\boldsymbol{\mu}+ {\bf O}{\bf Z};(n)}_{\pms}$the empirical distribution function of the sample $
      \boldsymbol{\mu}+ {\bf O}{\bf Z}_1,\ldots, \boldsymbol{\mu}+~\!{\bf O}{\bf Z}_n$ 
       associated with the grid ${\bf O}\mathfrak{G}_n$ (hence, by~${\bf F}^{{\bf Z};(n)}_{\pms}$ the empirical distribution function of the sample~${\bf Z}_1,\ldots, {\bf Z}_n$ associated with the grid $\mathfrak{G}_n$), 
      \end{compactenum} 
      \begin{equation}\label{equiv}
      {\bf F}^{\boldsymbol{\mu}+ {\bf O}{\bf Z};(n)}_{\pms} (\boldsymbol{\mu} + {\bf O}{\bf Z}_i)= {\bf O}{\bf F}^{{\bf Z};(n)}_{\pms}({\bf Z}_i), \quad  i=1,\ldots,n;
      \end{equation}
\begin{compactenum}
   \item[{\it (iii)}]  the center-outward ranks $R\n_{i;{\pms}}$ and the cosines ${\bf S}_{i;{\pms}}^{(n)\prime}{\bf S}_{j;{\pms}}\n$ computed from the sample~${\bf Z}_1,\ldots,{\bf Z}_n$ and the grid $\mathfrak{G}_n$ are the same as those computed from the \linebreak sample $\boldsymbol{\mu}+ {\bf O}{\bf Z}_1,\ldots, \boldsymbol{\mu}+ {\bf O}{\bf Z}_n$ and the grid ${\bf O}\mathfrak{G}_n$. 
   \end{compactenum}
\end{Proposition}

See Appendix~A.2 for the proof. 

These orthogonal equivariance and invariance  properties, however, do not extend to non-orthogonal affine transformations.  

\section{  H\' ajek representation and asymptotic normality}\label{Rpm} 
As in H\' ajek and \v{S}id\' ak~(1967), the rank-based statistics to be used in this context are quadratic forms in vectors of linear rank statistics---involving center-outward ranks and signs instead of ordinary ranks, though. Fundamental in H\' ajek's approach is an asymptotic representation result establishing the asymptotic equivalence between linear rank statistics and sums of independent variables. We start with a center-outward version of that result; asymptotic normality follows as a corollary. 

\subsection{Linear center-outward rank statistics} \label{sec.3.1}   
Linear rank statistics in this context depend on a {\it score function} ${\bf J}: {\mathbb{S}_d}\to\R^d$ and are indexed by  triangular arrays $\{c\n_1,\ldots,c\n_n\}$  
of real numbers (regression constants). On those score functions and regression constants we are making the following assumptions.  
\begin{Assumption}\label{Ass21}{\rm 
\begin{compactenum}
\item[{\it (i)}]  ${\bf J}: {\mathbb{S}_d}\to\R^d$ is continuous over ${\mathbb{S}_d}$;
\item[{\it \, (ii)}] for any sequence  ${\mathfrak s}\n=\{{\bf s}\n_1,\ldots,{\bf s}\n_n\}$ of $n$-tuples in ${\mathbb{S}_d}$ such that the uniform discrete distribution over ${\mathfrak s}\n$ converges weakly to ${\rm U}_d$ 
 as~$n\to~\!\infty$,
\begin{equation}\label{monass}\lim_{n\to\infty}n^{-1}
  \text{\rm tr} 
  \sum_{r=1}^n
  {\bf J}({\bf s}\n_r){\bf J}\pr({\bf s}\n_r)=
 \text{\rm tr} 
 \int_{\mathbb{S}_d}{\bf J}({\bf u}){\bf J}\pr({\bf u})\,
    {\mathrm d}{\rm U}_d 
    \end{equation}{where $ \int_{\mathbb{S}_d}{\bf J}({\bf u}){\bf J}\pr({\bf u})\,
    {\mathrm d}{\rm U}_d <\infty$ has full rank.} 
\end{compactenum}
}\end{Assumption}

As we shall see, a special role is played, in relation with spherical distributions,  by score functions of the form 
\begin{equation}\label{bfJ}
 {\mathbf J}({\bf u}):=J(\Vert {\bf u}\Vert)\frac{{\bf u}}{\Vert{\bf u}\Vert} {\bf 1}_{[\Vert{\bf u}\Vert\neq 0]} \qquad
  {\bf u}\in\mathbb{S}_d
\end{equation}
for some function $J: [0,1)\to \R$. Assumption~\ref{Ass21}  then holds if {\it (i)} $J$ is continuous and {\it (ii)}  
    \begin{equation}\label{Jsphcond} 0<{\lim_{n\to\infty}n^{-1}\sum_{r=1}^nJ^2\big(r/(n+1)\big)=\int_0^1J^2(u)\,
    {\mathrm d}u <\infty}
        \end{equation}
    (a sufficient condition for \eqref{Jsphcond} is the traditional assumption that $J$ has bounded variation,  {i.e.,}~is the difference of two nondecreasing functions).  Both \eqref{monass} and \eqref{Jsphcond} extend the conditions on univariate scores in Section~V.1.6 of \cite{HajekSidak}.  
    
As for the regression constants, we assume that the  classical {\it Noether conditions}   hold. 
 \begin{Assumption}\label{Ass22} {\rm The $c\n_i$'s are not all equal (for given $n$)  and
  satisfy 
\begin{equation}\label{Noeth}
\sum_{i=1}^n (c\n_i-\barc\n)^2 /{\max\limits_{1\leq i\leq n}(c\n_i-\barc\n)^2}\longrightarrow \infty\quad\text{ as } n\to\infty \vspace{-2mm}
\end{equation}
where~$
{ \barc\n:=n^{-1} \sum_{i=1}^n c\n_i}$. 
}\end{Assumption}

Associated with the score functions ${\bf J}$,  
 consider the $d$-dimensional  statistics \vspace{-1mm}
\begin{align} 
  \tenq{\tT}\n_a 
 &={\Big({\sum_{i=1}^n (c_i\n-\barc\n)^2}\Big)^{-1/2}}\sum_{i=1}^n (c_i\n-\barc\n) {\mathbf J}(\Fs(\tZ_i\n)),
\label{Tapp} 
\end{align}\vspace{-7mm}
\begin{align*}
   \tenq{\tT}\n_e &:=
    {\Big({\sum_{i=1}^n (c_i\n-\barc\n)^2}\Big)^{-1/2}}
    \sum_{i=1}^n
     (c_i\n-\barc\n) 
      \E\left[ {\mathbf J}({\bf F}_\pms(\tZ_i\n))  \bigg\vert \Fs(\tZ_i\n) \right], 
\nonumber\end{align*}
and\vspace{-4mm}
\begin{align*}
{\tT}\n &:=
  {\Big({\sum_{i=1}^n (c_i\n-\barc\n)^2}\Big)^{-1/2}}\sum_{i=1}^n (c_i\n-\barc\n) {\mathbf J}({\bf F}_\pms(\tZ_i\n)).
\end{align*}
Adopting H\' ajek's terminology, call $\tenq{\tT}\n_a $  an {\it approximate-score} linear rank statistic and $\tenq{\tT}\n_e$ an {\it exact-score} linear rank statistic. As we shall see, both~$\tenq{\tT}\n_a $ and $\tenq{\tT}\n_e$ admit the same asymptotic representation~${\tT}\n\!$, hence  are asymptotically equivalent. For score functions of the form~\eqref{bfJ}, we have
\[  \tenq{\tT}\n_a ={\Big({\sum_{i=1}^n (c_i\n-\barc\n)^2}\Big)^{-1/2}}\sum_{i=1}^n (c_i\n-\barc\n) J\Big(\frac{R\n_{i;{\pms}}}{n_R+1}\Big){\bf S}\n_{i;{\pms}}  
,\]
\begin{align*}  \tenq{\tT}\n_e =&{\Big({\sum_{i=1}^n (c_i\n-\barc\n)^2}\Big)^{-1/2}}
\\
  & \times
   \sum_{i=1}^n (c_i\n-\barc\n)
    \E\left[J\big(\big\Vert {\bf F}_\pms(\tZ_i\n)\big\Vert \big)
    \frac{ {\bf F}_\pms(\tZ_i\n) }{\big\Vert {\bf F}_\pms(\tZ_i\n)\big\Vert}
    \bigg\vert \Fs(\tZ_i\n) \right]
,\end{align*}
and
\[{\tT}\n ={\Big({\sum_{i=1}^n (c_i\n-\barc\n)^2}\Big)^{-1/2}} \sum_{i=1}^n (c_i\n-\barc\n) J(\big\Vert {\bf F}_\pms(\tZ_i\n)\big\Vert ) \frac{{\bf F}_\pms(\tZ_i\n)}{\big\Vert {\bf F}_\pms(\tZ_i\n)\big\Vert }\
.\]
\subsection{Asymptotic representation and asymptotic normality}

The following proposition is a center-outward multivariate counterpart    of the asymptotic results in Section~V.1.6 of \cite{HajekSidak}. Throughout this section, we assume that 
   ${\bf F}\n_\pms$ is computed from a triangular array $(\tZ_1\n,\dots, \tZ_n\n)$, $n\in\mathbb{N}$ of i.i.d.\ $d$-dimensional random vectors with distribution $\mathrm{P}\in{\mathcal P}_d^+$ and center-outward distribution function ${\bf F}_\pms$; the notation\   $o_{\text{\rm q.m.}}(1)$ is used for a sequence of random vectors
  tending to
zero in quadratic mean (hence also in probability). 

\begin{Proposition}[H\' ajek representation]\label{Prop22} Let Assumptions~\ref{Ass21} and \ref{Ass22} 
  hold \linebreak and~${\bf Z}\n_1,\ldots,{\bf Z}\n_n$ be  i.i.d.~with distribution ${\rm P}\in{\cal P}^{+}_{d}$. Then, 
 \[{\it (i)}\ \tenq{\tT}\n_a - {\tT}\n =o_{\text{\rm q.m.}}(1),\quad {\it (ii)}\ \tenq{\tT}\n_e - {\tT}\n =o_{\text{\rm q.m.}}(1),
\quad\text{and}\quad 
{\it (iii)}\  \tenq{\tT}\n_a - \tenq{\tT}\n_e =o_{\text{\rm q.m.}}(1)
\]
   as~$n\to\infty$ in such a way that $n_R\to\infty$ and
  $n_S\to\infty$.
\end{Proposition}

See Appendix~A.3 for the proof.

 The asymptotic normality of $\tenq{\tT}\n_a$ and $\tenq{\tT}\n_e$  
then follows from {Proposition} \ref{Prop22} and  the asymptotic normality of ${\tT}\n $, along with the distribution-freeness of~$\tenq{\tT}\n_a$ and $\tenq{\tT}\n_e$.

\begin{Proposition}[Asymptotic normality]\label{asN} Let Assumptions~\ref{Ass21} and~\ref{Ass22}
  hold and ${\bf Z}\n_1,\ldots,{\bf Z}\n_n$ be  i.i.d.~with distribution ${\rm P}\in{\cal P}_d$. Then, 
  $\tenq{\tT}\n_a$, 
 $\tenq{\tT}\n_e\vspace{1mm}$, and
  ${\tT}\n $ are asymptotically normal  as~$n\to~\!\infty$ (in such a way that~$n_R\to\infty$ and~$n_S\to\infty$), with mean $\mathbf 0$  and
  covariance~$\int_{\mathbb{S}_d}{\bf J}({\bf u}){\bf J}\pr({\bf u})\,
    {\mathrm d}{\rm U}_d$ reducing, for~$\bf J$ of the form~\eqref{bfJ}, to  
    $d^{-1}\int_0^1J^2(u)\, \mathrm{d}u\, {\mathbf I}_d$.   
\end{Proposition}

See Appendix~A.4 for the proof.

\setcounter{equation}{0}

\section{Multiple-output  linear models}
 Based on the center-outward ranks and signs  of Section~\ref{Rpm}, we now construct rank  tests for the slopes of multiple-output linear models, extending to a multivariate setting the methods developed,  {e.g.,}~in \cite{PuriSen1985} for the single-output case. 

\subsection{The model}
Consider the multiple-output linear (or multiple-output regression) model  under which an observed ${\mathbf Y}\n$ satisfies  
\begin{equation}\label{MRegr}{\mathbf Y}\n=
  {\mathbf 1}_n{\boldsymbol\beta}_0\pr +
 {\mathbf C}\n{\boldsymbol\beta} + {\boldsymbol\varepsilon}\n,
\end{equation}
where $  {\mathbf 1}_n:=(1,\ldots,1)\pr$,  
\[
{\mathbf Y}\n=\left(\begin{array}{cccc}
 Y\n_{11}& Y\n_{12}&\ldots&Y\n_{1d}\\ 
\vdots &\vdots&&\vdots \\ 
Y\n_{n1}& Y\n_{n2}&\ldots&Y\n_{nd}
\end{array}\right)
=
\left(\begin{array}{c} {\mathbf Y}^{(n)\prime}_1\\ \vdots \\ {\mathbf Y}^{(n)\prime}_n
\end{array}
\right)
\]
is an $n\times d$ matrix of $n$ observed $d$-dimensional outputs, 
\[
{\mathbf C}\n=\left(\begin{array}{cccc}
 c\n_{11}& c\n_{12}&\ldots&c\n_{1m}\\ 
\vdots &\vdots&&\vdots \\ 
c\n_{n1}& c\n_{n2}&\ldots&c\n_{nm}
\end{array}\right)
=
\left(\begin{array}{c} {\mathbf c}^{(n)\prime}_1\\ \vdots \\ {\mathbf c}^{(n)\prime}_n
\end{array}
\right)
\]
 an $n\times m$ matrix of (specified) deterministic covariates,  
 \[
 \boldsymbol\beta_0\pr=(\beta_{01},\ldots,\beta_{0d})
 \quad
 \text{and}
 \quad
 {\boldsymbol\beta}=
\left(\begin{array}{cccc}
 \beta_{11}& \beta_{12}&\ldots&\beta_{1d}\\ 
\vdots &\vdots&&\vdots \\ 
\beta_{m1}& \beta_{m2}&\ldots&\beta_{md}
\end{array}\right)
=
\left(\begin{array}{c} {\boldsymbol\beta}^{\prime}_1\\ \vdots \\ {\boldsymbol\beta}^{\prime}_m
\end{array}
\right)
 \]
  a $d$-dimensional intercept and   an $m\times d$ matrix of regression coefficients, and 
\[
{\boldsymbol\varepsilon}\n=\left(\begin{array}{cccc}
 \varepsilon\n_{11}& \varepsilon\n_{12}&\ldots&\varepsilon\n_{1d}\\ 
\vdots &\vdots&&\vdots \\ 
\varepsilon\n_{n1}& \varepsilon\n_{n2}&\ldots&\varepsilon\n_{nd}
\end{array}\right)
=
\left(\begin{array}{c} {\boldsymbol\varepsilon}^{(n)\prime}_1\\ \vdots \\ {\boldsymbol\varepsilon}^{(n)\prime}_n
\end{array}
\right)
\]
an $n\times d$ matrix of nonobserved i.i.d.~$d$-dimensional errors ${\boldsymbol\varepsilon}^{(n)}_i$, $i=1,\ldots,n$ with density~$f^{\boldsymbol\varepsilon}$.  If 
$ \boldsymbol\beta_0$ is to be identified, a location constraint has to be imposed on $f^{\boldsymbol\varepsilon}$. One could think of   the classical constraint ${\rm E}{\boldsymbol\varepsilon}^{(n)}_i ~\!=~\!{\bf 0}$ (requiring the existence   of a finite mean): $ \boldsymbol\beta_0 +  \boldsymbol\beta\pr{\bf c}\n_i$ then is to be interpreted as the expected value of ${\bf Y}\n_i$ for covariate values ${\bf c}\n_i$. In the context of this paper, however,  a more natural location constraint (which moreover does not require any integrability condition) is   ${\bf F}^{{\boldsymbol\varepsilon}}_{\pms}({\bf 0}) ={\bf 0}$, where ${\bf F}^{{\boldsymbol\varepsilon}}_{\pms}$ stands for the center-outward distribution function of the ${\boldsymbol\varepsilon}^{(n)}_i$'s: $\bf 0$ and $ \boldsymbol\beta_0 +  \boldsymbol\beta\pr{\bf c}\n_i$  then are  {\it center-outward medians} for $\boldsymbol \varepsilon$ and~${\bf Y}\n_i$, respectively.  
 
  In most applications, however, one is interested mainly in the impact of the input covariates~${\bf c}\n_i$ on the output~${\bf Y}\n_i$:  the matrix $ \boldsymbol\beta$  is the parameter of interest, and $ \boldsymbol\beta_0$ is a nuisance. There is no need, then, for identifying~$\boldsymbol\beta_0$ nor qualifying $ \boldsymbol\beta_0 +  \boldsymbol\beta\pr{\bf c}\n_i$  as a  mean or a center-outward median for ${\bf Y}\n_i$:~$\boldsymbol \beta$ is to be interpreted as a matrix of treatment effects governing the shift~${\boldsymbol\delta}\pr{\boldsymbol\beta}$ in the distribution of the $d$-dimensional output produced by a variation ${\boldsymbol\delta}$ in the~$m$-dimensional covariate.  Center-outward ranks and signs being insensitive to shifts, there is even no need to specify, nor to estimate $ \boldsymbol\beta_0$. 
  
  
  \subsection{Local Asymptotic Normality (LAN)} The model \eqref{MRegr} is easily seen to be locally asymptotically normal (LAN) under the following two classical assumptions.\vspace{-1mm} 
   
  \begin{Assumption}\label{Ass30} {\rm The square root ${\bf z}\mapsto\left(f^{\boldsymbol\varepsilon}\right)^{1/2}({\bf z})$ of the error density is  {\it differentiable in quadratic mean},\footnote{It follows from a result\vspace{-2mm}
   by \cite{Rouss72} independently rediscovered by \cite{GH95} that quadratic mean differentiability   is equivalent to  partial quadratic mean derivability\vspace{-2mm} with respect to all variables.} with quadratic mean gradient $\nabla\left(f^{\boldsymbol\varepsilon}\right)^{1/2}\!$. Letting $
  {\boldsymbol\varphi}_{f^{\boldsymbol\varepsilon}}\!:=-2\nabla\left(f^{\boldsymbol\varepsilon}\right)^{1/2}\!/\!\left(f^{\boldsymbol\varepsilon}\right)^{1/2}
\!\!$, 
   assume moreover that the {\em information matrix} ${\boldsymbol{\mathcal I}}_{f^{\boldsymbol\varepsilon}}\!:={\rm E}\left[{\boldsymbol\varphi}_{f^{\boldsymbol\varepsilon}}({\boldsymbol\varepsilon}){\boldsymbol\varphi}\pr_{f^{\boldsymbol\varepsilon}}({\boldsymbol\varepsilon})\right] $ has full rank $d$.
}  \end{Assumption}

On the regression constants ${\mathbf C}\n$,  we borrow from \cite{HallinPain2005}  
 the following assumptions; note that Part~{\it (iii)} requires that each of the $m$ triangular arrays of constants~$c\n_{ij}$, $i\in\mathbb{N}$, $j=1,\ldots, m$ satisfies Assumption~\ref{Ass22}.

\begin{Assumption}\label{Ass31}{\rm Let~$\bar{\bf c}\n\!:=n^{-1}\!\sum_{i=1}^n{\bf c}\n_i\!$,  ${\bf V}_{\bf c}\n\!:=n^{-1}\sum_{i=1}^n\big({\bf c}\n_i\! - \bar{\bf c}\n\big)\big({\bf c}^{(n)}_i\!  - \bar{\bf c}\n\big)^\prime$, and denote   
  by ${\bf D}_{\bf c}\n$ the diagonal matrix with diagonal  elements~$\big({\bf V}_{\bf c}\n\big)_{jj}$,  $j=1,\ldots,m$:
\begin{compactenum}
\item[{\it (i)}]  $\big({\bf V}_{\bf c}\n\big)_{jj}>0$ for   $j=1,\ldots,m$;
\item[{\it (ii)}] defining ${\bf R}_{\bf c}\n:={\bf D}_{\bf c}^{(n)-1/2}{\bf V}_{\bf c}\n{\bf D}_{\bf c}^{(n)-1/2}$, the limit ${\bf R}_{\bf c}:=\lim_{n\to\infty}{\bf R}_{\bf c}\n$ exists, is positive definite, and factorizes into ${\bf R}_{\bf c}=\big({\bf K}_{\bf c}{\bf K}_{\bf c}^{\prime}\big)^{-1}$ for some full-rank $m\times m$ matrix ${\bf K}_{\bf c}$;
\item[{\it (iii)}] letting $\bar{c}\n_j:=n^{-1}\sum_{i=1}^nc_{ij}\n$,   the following Noether conditions hold: \vspace{-2mm}
\[ \lim_{n\to\infty}\sum_{i=1}^n \big(c\n_{ij} -\bar{c}\n_j \big)^2/\max_{1\leq i\leq n}\big(c\n_{ij} -\bar{c}\n_j \big)^2=\infty ,
\qquad \text{$j=1,\ldots,m$}.\]
\end{compactenum}
}\end{Assumption}

Letting ${\bf Z}\n_i= {\bf Z}\n_i({\boldsymbol\beta}):={\bf Y}\n_i 
\!\!-  {\mathbf 1}_n{\boldsymbol\beta}_0\pr  - {\boldsymbol\beta}^\prime{\bf c}\n_i$, the following result  readily follows from, e.g., \cite[Theorem~12.2.3]{LR05}. In order to simplify the notation, we throughout adopt the same contiguity rates as in \cite{HallinPain2005}. Namely, we consider local perturbations of the parameter $\boldsymbol\beta$  of the form ${\boldsymbol\beta}+{\boldsymbol\nu}(n){\boldsymbol\tau}$ where $\boldsymbol\tau$ is an~$m\times d$   matrix and~${\boldsymbol\nu}(n):=n^{-1/2}{\bf K}\n_{\bf c}$, 
 with ${\bf K}\n_{\bf c}:= \big({\bf D}\n_{\bf c}\big)^{-1/2}{\bf K}_{\bf c}$. This, which incorporates the asymptotic behavior of the regression constants, is a notational convenience and has no impact on the form of locally asymptotically optimal test statistics.

\begin{Proposition}\label{LAN} Under Assumptions~\ref{Ass30} and \ref{Ass31}, the model \eqref{MRegr}  is LAN (with respect to ${\boldsymbol\beta}$), with central sequence $\boldsymbol\Delta^{(n)}_{{\boldsymbol\beta}_0 ;{f}^{\boldsymbol\varepsilon}}({\boldsymbol\beta}):=n^{1/2}\text{\rm vec}{\boldsymbol\Lambda}\n_{{\boldsymbol\beta}_0 ;{f}^{\boldsymbol\varepsilon}}$ where 
\begin{equation}\label{Deltagen} 
 {\boldsymbol\Lambda}\n_{{\boldsymbol\beta}_0 ;{f}^{\boldsymbol\varepsilon}}:=\frac{1}{n} \sum_{i=1}^n 
 {\bf K}\npr_{\bf c} \big({\bf c}\n_i - \bar{\bf c}\n
\big) {\boldsymbol\varphi}_{f^{{\boldsymbol\varepsilon}}}\pr({\bf Z}\n_i )
\end{equation}
and Fisher information ${\boldsymbol{\mathcal I}}_{f^{\boldsymbol\varepsilon}}\otimes{\bf I}_m$.
 \end{Proposition}

LAN for the same linear model \eqref{MRegr} has been established (in the broader context of regression with VARMA errors   
 in \cite{HallinPain2005})  
under the assumption that the error density $f^{{\boldsymbol\varepsilon}}$   is {\it centered elliptical} (for simplicity, we henceforth are dropping the word ``centered''), that is, has the form
\begin{equation}\label{ellf}
f^{{\boldsymbol\varepsilon}}({\bf z}) = \kappa^{-1}_{d,{\mathfrak f}}\big(\text{det}{\boldsymbol \Sigma}\big)^{-1/2}{\mathfrak f}\big( ({\bf z}^\prime{\boldsymbol \Sigma}^{-1}{\bf z})^{1/2}
\big)
\end{equation} 
with   $\kappa_{d,{\mathfrak f}}:= \big({2\pi^{d/2}}/{\Gamma (d/2)}\big)\int_0^\infty r^{d-1}{\mathfrak f}(r)\,{\mathrm d}r$ 
for some symmetric positive definite {\it shape matrix $\boldsymbol \Sigma$} and some {\it radial density}~$\mathfrak f$ 
 (over $\mathbb{R}^+_0$) such that~$\mathfrak f(z)>0$ Lebesgue-a.e.\ in~$\mathbb{R}^+_0$ and~$\int_0^\infty r^{d-1}{\mathfrak f}(r)\,{\mathrm d}r<~\!\infty$. When $\boldsymbol\varepsilon$ is elliptical with  shape matrix~$\boldsymbol \Sigma$ and radial density~$\mathfrak f$, 
 the  modulus
 $\Vert\boldsymbol \Sigma^{-1/2}\boldsymbol\varepsilon \Vert$ has probability density 
  ${\mathfrak f}^\star_{ d}(r)=(\mu_{d-1;{\mathfrak{f}}})^{-1}r^{d-1}\mathfrak{f}(r)I[r>0]$,  
  where~$\mu_{d-1;{\mathfrak{f}}}\!:=\!\int_0^\infty r^{d-1}\mathfrak{f}(r){\rm d}\, r$, 
 and   distribution function $F^\star_{d;{\mathfrak f}}$.
 
 Assumption~\ref{Ass30} then is equivalent to the  mean square differentiability, with quadratic mean derivative $\big({\mathfrak f}^{1/2}\big)\pr\!$,  of $x\mapsto{\mathfrak f}^{1/2}(x)$, $x\in\R_0^+$ (a scalar); 
 letting $\varphi_{\mathfrak f}\!:=\!-2\big({\mathfrak f}^{1/2}\big)\pr\!/ {\mathfrak f}^{1/2}\!
 $, we automatically get
  $
{\cal I}_{d;{\mathfrak f}}:=
 \int_0^1\!\Big(\varphi_{\mathfrak f} \circ \big( F^\star_{d;{\mathfrak f}}\big)^{-1}\!(u) \Big) ^2{\rm d}u <~\!\infty.$ 
  Define the {\it sphericized residuals} 
\begin{equation}
 {\bf Z}_i^{(n)\,\text{\rm ell}}\!\!
 :=\!\big(\widehat{\boldsymbol\Sigma}^{(n)}\big)^{-1/2}\big({\bf Y}\n_i\!\! -  {\boldsymbol\beta}_0  - {\boldsymbol\beta}^\prime{\bf c}\n_i  \big)
 =\big(\widehat{\boldsymbol\Sigma}^{(n)}\big)^{-1/2}\big({\bf Z}\n_i\big)
 ,\ \  i=1,\ldots,n
\label{Zell}\end{equation}
where   the matrix $\big(\widehat{\boldsymbol\Sigma}^{(n)} \big)^{1/2}$  is the symmetric root of a consistent estimator~$\widehat{\boldsymbol\Sigma}^{(n)}$ of some mul\-tiple $a\,{\boldsymbol\Sigma}$ of ${\boldsymbol\Sigma}$ ($a>0$ an arbitrary constant) satisfying the following consistency assumption.
\begin{Assumption}\label{Ass32} {\rm Under \eqref{MRegr}, $\widehat{\boldsymbol\Sigma}^{(n)}\!-a {\boldsymbol\Sigma}= O_{\rm P}(n^{-1/2})$ as $n\to\infty$, for some $a>0$; moreover, $\widehat{\boldsymbol\Sigma}^{(n)}$ is invariant under permutations and reflections (with respect to the origin) of the residuals  ${\bf Z}\n_i=({\bf Y}\n_i\! -  {\mathbf 1}_n{\boldsymbol\beta}_0\pr - {\boldsymbol\beta}^\prime{\bf c}\n_i)$'s, and equivariant under their affine transformations.
}\end{Assumption}

A traditional choice which, however, rules out heavy-tailed radial densities with infinite second-order moments, is the empirical covariance 
 of the~${\bf Z}\n_i$'s. An alternative,   satisfying Assumption~\ref{Ass32} without  any moment assumptions, is Tyler's estimator of scatter,  see Theorems~4.1 and 4.2  in~\cite{Tyl87} for its  strong consistency and asymptotic normality.   

Under Assumption~\ref{Ass32}, which entails the affine invariance of $ {\bf Z}_i^{(n)\,\text{\rm ell}}$, Proposition~\ref{LAN}  takes the following form.
 \begin{Proposition}\label{LANell} Under Assumptions~\ref{Ass31} and~\ref{Ass32}, the model \eqref{MRegr} with error density~$f^{\boldsymbol\varepsilon}$ of the elliptical  type~\eqref{ellf} and quadratic mean differentiable~${\mathfrak f}^{1/2}$ is LAN (with respect to~${\boldsymbol\beta}$), with central sequence\vspace{-1mm} 
\begin{equation}
\label{Deltaell} 
\boldsymbol\Delta^{(n)\,\text{\rm ell}}_{{\widehat{\boldsymbol\Sigma}\n}\! , {\boldsymbol\beta}_0 ;{\mathfrak f}}({\boldsymbol\beta}):=
 n^{1/2} \left( \big({\widehat{\boldsymbol\Sigma}\n}\big) ^{-1/2}\otimes{\bf I}_m\right)\text{\rm vec}{\boldsymbol\Lambda}^{(n)\,\text{\rm ell}}_{{\widehat{\boldsymbol\Sigma}\n}\! , {\boldsymbol\beta}_0 ;{\mathfrak f}}({\boldsymbol\beta})
=
 \boldsymbol\Delta^{(n)\,\text{\rm ell}}_{{{\boldsymbol\Sigma}} , {\boldsymbol\beta}_0 ;{\mathfrak f}}({\boldsymbol\beta})
+o_{\rm P}(1)
\end{equation}
 where\vspace{-3mm}  
\begin{align}\label{DeltaellLambd} 
{\boldsymbol\Lambda}^{(n)\,\text{\rm ell}}_{{\widehat{\boldsymbol\Sigma}\n}\! , {\boldsymbol\beta}_0 ;{\mathfrak f}}({\boldsymbol\beta}):=
\frac{1}{n}\sum_{i=1}^n
\varphi_{\mathfrak f}\big(\big\Vert  {\bf Z}_i^{(n)\,\text{\rm ell}}   \big\Vert\big)
 {\bf K}\npr_{\bf c} \big({\bf c}\n_i - \bar{\bf c}\n
\big) \left(\frac{ {\bf Z}_i^{(n)\,\text{\rm ell}} }{\big\Vert  {\bf Z}_i^{(n)\,\text{\rm ell}}   \big\Vert} \right)\pr ,
\end{align}

\noindent yielding  a Fisher information matrix  $\frac{1}{d}{\cal I}_{d;{\mathfrak f}}\,{\boldsymbol\Sigma}^{-1} \otimes {\bf I}_{m}$.
\end{Proposition}
This LAN result, where the residuals are subjected to  preliminary (empirical) sphericization via $\big({\widehat{\boldsymbol\Sigma}\n}\big) ^{-1/2}$, stresses the fact that elliptical families with given~$\mathfrak f$  are parametrized  {\it spherical families} (indexed by 
 $\boldsymbol\Sigma$). Actually, since 
 $\boldsymbol\Delta^{(n)\,\text{\rm ell}}_{{{\boldsymbol\Sigma}} , {\boldsymbol\beta}_0 ;{\mathfrak f}}({\boldsymbol\beta})=\left({\bf I}_d\otimes {{\boldsymbol\Sigma}} ^{-1/2}\right)\boldsymbol\Delta^{(n)\,\text{\rm ell}}_{{{\bf I}_d} , {\boldsymbol\beta}_0 ;{\mathfrak f}}({\boldsymbol\beta}),$ 
 the limiting Gaussian shift experiments associated with elliptical and spherical errors coincide (with the perturbation $\text{\rm vec}({\boldsymbol\tau})$ of $\text{\rm vec}({\boldsymbol\beta})$ in the elliptical case corresponding to a perturbation $\text{\rm vec}({\boldsymbol\varsigma})=\left({\bf I}_d\otimes {{\boldsymbol\Sigma}} ^{-1/2}\right)\text{\rm vec}({\boldsymbol\tau})$ in the spherical case).  That invariance under linear sphericization  of local limiting Gaussian shifts, however, does not extend to the general case of Proposition~\ref{LAN}.

  \color{black}
  \setcounter{equation}{0}

  \section{Rank tests for multiple-output linear models}\label{Rpmsec}
  \subsection{Elliptical (Maha\-la\-nobis)  rank tests }\label{sec_5.1}

Rank-based inference for elliptical multiple-output linear models was developed in \cite{HallinPain2005}.  The ranks  and the signs  there are the {\it elliptical} or~{\it Maha\-la\-nobis  ranks  and signs}---namely, 
 the ranks  $R_i^{(n)\,\text{\rm{ell}}}$  of the moduli $\big\Vert {\bf Z}_i^{(n)\,\text{\rm{ell}}}\Vert$ and  the signs (directions)~${\bf S}_i^{(n)\,\text{\rm{ell}}\!\!}:= {\bf Z}_i^{(n)\,\text{\rm{ell}}}\!/\Vert{\bf Z}_i^{(n)\,\text{\rm{ell}}}\Vert$, both computed, in agreement with the above remark on the spherical nature of elliptical families,  after the empirical sphericization~\eqref{Zell}.

\textcolor{black}{
Consider the null hypothesis $H\n_0 ({\boldsymbol\beta}^0)$ under which ${\bf Y}\n$ satisfies \eqref{MRegr} with ${\boldsymbol\beta}={\boldsymbol\beta}^0$, specified ${\boldsymbol\beta}_0$, elliptical $f^{{\boldsymbol\varepsilon}}$, and radial density ${\mathfrak f}$. \cite{HallinPain2005} define
\[
 \tenq{\boldsymbol\Lambda}^{(n)\,\text{\rm{ell}}}_J:=n^{-1}
\sum_{i=1}^nJ\Big(\frac{R_i^{(n)\,\text{\rm{ell}}}}{n+1}\Big)\, {\bf S}_i^{(n)\,\text{\rm{ell}}}\big({\bf c}\n_i - \bar{\bf c}\n
\big)\pr {\bf K}\n_{\bf c}
\]
for a score function $J:[0,1)\to\mathbb{R}$ and show that the test of $H\n_0 ({\boldsymbol\beta}^0)$  can be based on  
\[
\utilde{Q}^{(n)\,\text{\rm ell}}_{J}({\boldsymbol\beta}^0)=\frac{n\, d}{\int_0^1 J^2(u)du} \big(\text{vec} \utilde{{\boldsymbol\Lambda}}^{(n)\,\text{\rm ell}}_{J}\big)\pr\big(\text{vec} \utilde{{\boldsymbol\Lambda}}^{(n)\,\text{\rm ell}}_{J}\big)\label{utQ},
\]
which is asymptotically chi-square with $md$ degrees of freedom under the null. 
 }

The validity of tests based on those elliptical ranks and signs, unfortunately, requires an elliptical $f^{{\boldsymbol\varepsilon}}$. 
 A welcome relaxation of stricter Gaussianity assumptions, ellipticity remains an extremely  strong symmetry requirement; it is made, essentially, for lack of anything better but is unlikely to hold in practice.  If  the assumption of ellipticity is to be waived, elliptical ranks and signs are losing their distribution-freeness   for the benefit of the center-outward ranks and signs. And, since  center-outward ranks and signs, in view of Proposition~\ref{invF}, are  invariant under location shift,   center-outward rank tests can address the (more realistic) unspecified intercept case without any additional estimation step. 

  \subsection{Center-outward  rank tests }
Denote by ${\bf F}\n_\pms$ the empirical center-outward distribution associated with the observed $n$-tuple $({\bf Z}\n_1,\ldots, {\bf Z}\n_n)$ where   ${\bf Z}\n_i$ now is defined as ${\bf Y}\n_i\! -  {\boldsymbol\beta}^\prime{\bf c}\n_i$,  by~$R\n_{i;{\pms}}$ and ${\bf S}\n_{i;{\pms}}$, respectively, the corresponding center-outward ranks and signs. In line with the form of the central sequence  \eqref{Deltagen}, 
 consider \vspace{-1mm}
\begin{equation} \label{tildeLambda}
\tenq{\boldsymbol\Lambda}^{(n)\pms}_{\mathbf J}:=n^{-1}
\sum_{i=1}^n {\bf K}\npr_{\bf c} \big({\bf c}\n_i - \bar{\bf c}\n
\big) {\mathbf J}\pr\Bigg(\frac{R^{(n)}_{i;{\pms}}}{n_R+1}\, {\bf S}_{i;{\pms}}^{(n)}\Bigg). \vspace{-1mm}
\end{equation}

It follows from the asymptotic representation result of
Proposition~\ref{Prop22} that, when the actual density is
$f^{\boldsymbol\varepsilon}$, for the scores ${\bf J}=
{\boldsymbol\varphi}_{f^{{\boldsymbol\varepsilon}}}\circ{\bf
  F}_\pm^{-1}$, with
${\boldsymbol\varphi}_{f^{{\boldsymbol\varepsilon}}}$ defined in Assumption~\ref{Ass30}   \vspace{-1mm}
\begin{equation} \label{tildeLambdaDelta}
\tenq{\boldsymbol\Delta}^{(n)}_{{\boldsymbol\beta}_0 ;{f}^{\boldsymbol\varepsilon}}({\boldsymbol\beta}):=n^{1/2}\text{\rm vec}\tenq{\boldsymbol\Lambda}^{(n)\pms}_{\bf J} = \boldsymbol\Delta^{(n)}_{{\boldsymbol\beta}_0 ;{f}^{\boldsymbol\varepsilon}}({\boldsymbol\beta})+o_{\rm P}(1) \vspace{-1mm}
\end{equation}
 and $\tenq{\boldsymbol\Delta}^{(n)}_{{\boldsymbol\beta}_0 ;{f}^{\boldsymbol\varepsilon}}({\boldsymbol\beta})$ thus constitutes a version, based on the center-outward ranks and signs and hence distribution-free, of the central sequence  $\boldsymbol\Delta^{(n)}_{{\boldsymbol\beta}_0 ;{f}^{\boldsymbol\varepsilon}}({\boldsymbol\beta})$ in~\eqref{Deltagen}. 

{Recall that~$H\n_0 ({\boldsymbol\beta}^0)$ denotes the null hypothesis under which~${\boldsymbol\beta}={\boldsymbol\beta}^0$ while~${\boldsymbol\beta}_0$ and the distribution
${\rm P}\in {\cal P}_d$
of the~$\boldsymbol\varepsilon$'s remains unspecified; denote by $H\n_1({\boldsymbol\beta}^0, {\bf B},f)$ any local sequence of alternatives under which ${\boldsymbol\beta}={\boldsymbol\beta}^0+{\boldsymbol\nu}(n){\bf B}\n={\boldsymbol\beta}^0+ n^{-1/2}{\bf K}_{\bf c}\n{\bf B}\n$ with~$n^{-1/2}{\bf K}\n_{\bf c}\left( {\bf B}\n -{\bf B}\right)=o(1)$ for some ${\bf B}\neq{\bf 0}$ (which entails contiguity)  and the errors~$\boldsymbol\varepsilon$ have density $f$.  Also recall   that a sequence of tests in a LAN experiment is called  {\it locally  asymptotically maximin at asymptotic level~$\alpha$}  if its power function converges pointwise, as $n\to\infty$, to the power function of an $\alpha$-level maximin test in the corresponding limit Gaussian shift experiment: see, e.g., Section~11.9 in \cite{LeCam86}. The following asymptotic   results then hold. 
}
\begin{Proposition}\label{asNpm} 
{ Let ${\bf Y}\n_i$ satisfy \eqref{MRegr}. Then, under  
 Assumptions~\ref{Ass21} and~\ref{Ass31}, 
\begin{compactenum}
\item[{\it (i)}]   
  $n^{1/2}\text{\rm vec}\tenq{\boldsymbol\Lambda}^{(n)\pms}_{\bf J}$ is    asymptotically normal, with mean~$\bf 0$ under the null hypothesis
  ~$H\n_0 ({\boldsymbol\beta}^0)$, mean 
  \begin{equation}\label{shifteq}
  {\boldsymbol\mu}_{f,{\bf B}}=
 \left(\int_{{\mathbb S}_d} {\bf J}({\bf u}) {\boldsymbol\varphi}_{f}'\left({\bf Q}_\pms(\bf{u})\right)\, {\rm d}{\rm U}_d ({\bf u}) \otimes {\bf I}_m \right) \text{\rm vec}\,{\bf B}
  \end{equation}
  under local 
  alternatives of the form $H\n_1({\boldsymbol\beta}^0, {\bf B},f)$, 
    and covariance
\begin{equation}\label{coveq}
{\boldsymbol{\mathcal I}}_{{\bf J}}\otimes{\bf I}_m:=  
\int_{{\mathbb S}_d}{\bf J}({\bf u}){\bf J}\pr({\bf u})\,{\rm d}{\rm U}_d  ({\bf{u}}) \otimes{\bf I}_m
\vspace{-3mm}\end{equation}
 under both;
\item[{\it (ii)}]the test rejecting  $H\n_0({\boldsymbol\beta}^0)$ whenever the test statistic
\begin{equation}\label{Qpm}
\tenq{{Q}}^{(n)\pms}_{\bf J}:=
n\left(\text{\rm vec}\tenq{\boldsymbol\Lambda}^{(n)\pms}_{\bf J}
\right)\pr
\left({\boldsymbol{\mathcal I}}^{-1}_{{\bf J}}\otimes{\bf I}_m\right)
 \left(\text{\rm vec}\tenq{\boldsymbol\Lambda}^{(n)\pms}_{\bf J}
\right)
\end{equation}
exceeds the $(1-\alpha)$ quantile $\chi^2_{md;1-\alpha}$ of a chi-square distribution with $md$ degrees of freedom has asymptotic level $\alpha$ as $n\to\infty$;\footnote{Since\vspace{-0.5mm} $\tenq{{\bf Q}}^{(n)\pms}_{\bf J}$ is distribution-free under the null hypothesis $H\n_0 ({\boldsymbol\beta}^0)\vspace{-1mm}$,  the finite-$n$ size of this test is uniform over $H\n_0 ({\boldsymbol\beta}^0)$, hence
  uniformly close to $\alpha$ for $n$ large enough.
  This is in sharp\vspace{-2mm} contrast with daily practice pseudo-Gaussian tests, which remain asymptotically valid under a broad range\vspace{-2mm} of distributions, albeit not uniformly so (see Section~\ref{pseudoSec}).  
    } its asymptotic power against alternatives of the form  $H\n_1({\boldsymbol\beta}^0, {\bf B},f)$ is 
    $1-F_{\chi^2_{md;q}} (\chi^2_{md;1-\alpha})
    $ 
    where $F_{\chi^2_{md;q}}$ stands for the noncentral chi-square distribution function with $md$ degrees of freedom and noncentrality parameter 
\[  {
    q\!=\! \text{\rm vec}^\prime {\bf B}\left[\!\left(\int_{{\mathbb S}_d}\!\!  {\boldsymbol\varphi}_{f}\!\left({\bf Q}_\pms(\bf{u})\right) {\bf J}' ({\bf u})\,  {\rm d}{\rm U}_d ({\bf{u}})\right)\! {\boldsymbol{\mathcal I}}_{{\bf J}}^{-1} \!\left(\int_{{\mathbb S}_d} \!\!{\bf J} ({\bf u})  {\boldsymbol\varphi}_{f}'\!\left({\bf Q}_\pms({\bf{u}})\right)\,  {\rm d}{\rm U}_d ({\bf{u}}) \right)\! \otimes  {\bf I}_m\!\right]\! \text{\rm vec\,}{\bf B};
 }   \]
  \item[{\it (iii)}] for ${\bf J}={\bf J}_{f^{\boldsymbol\varepsilon}}
  :={\boldsymbol\varphi}_{f^{{\boldsymbol\varepsilon}}}\circ{\bf F}_{{\boldsymbol\varepsilon}\pm}^{-1}$ where ${\bf F}_{{\boldsymbol\varepsilon}\pm}$ denotes the center-outward distribution function associated with $f^{\boldsymbol\varepsilon}$,  the covariance \eqref{coveq} coincides with~${\boldsymbol{\mathcal I}}_{f^{\boldsymbol\varepsilon}}\otimes{\bf I}_m$  
 and the test based on~$\tenq{{Q}}^{(n)\pms}_{{\bf J}_{f^{\boldsymbol\varepsilon}}}$ 
 (as described in {\it (ii)}) is   locally  asymptotically maximin,  at asymptotic level~$\alpha$, for   $H\n_0({\boldsymbol\beta}^0)$ against~$H\n_1({\boldsymbol\beta}^0, {\bf B},f^{\boldsymbol\varepsilon})$ at asymptotic level $\alpha$. 
\end{compactenum}
}\end{Proposition}

\begin{Corollary}\label{scorell}
\begin{compactenum}
\item[{\it (i)}] In the particular case of a {\em spherical score} of  
 the  form~\eqref{bfJ}, the test statistic $\tenq{{Q}}^{(n)\pms}_{\bf J}$  simplifies into\vspace{-2mm} 
\begin{equation}\label{Qpmell}
\tenq{{Q}}^{(n)\pms}_{J}=\frac{nd}{\int_{0}^1J^2(u){\rm d}u}
\text{\rm vec}^\prime\tenq{\boldsymbol\Lambda}^{(n)\pms}_{J}
 \text{\rm vec}\tenq{\boldsymbol\Lambda}^{(n)\pms}_{J}
\end{equation}
where 
 $\tenq{\boldsymbol\Lambda}^{(n)\pms}_{J}:= n^{-1}
\sum_{i=1}^n J\left(\frac{R\n_{i;\pm}}{n_R+1}
\right){\bf K}_{\bf c}\npr \big({\bf c}\n_i - \bar{\bf c}\n
\big){\bf S}_{i;{\pms}}\npr $ 
 and~$n^{1/2}\text{\rm vec}\tenq{\boldsymbol\Lambda}^{(n)\pms}_{J}$under~$H\n_0 ({\boldsymbol\beta}^0)$  is asymptotically normal with mean $\bf 0$ and vari\-ance~$d^{-1}{\int_{0}^1J^2(u){\rm d}u}\,{\bf I}_{md}$. 
\item[{\it (ii)}] The test statistic $\tenq{{Q}}^{(n)\pms}_{J_{\mathfrak f}}$ with spherical score 
$J_{\mathfrak f}:=\varphi_{\mathfrak f}\circ \big( F^\star_{d;{\mathfrak f}}\big)^{-1}
$ 
yields locally asymptotically optimal tests under the spherical density with radial density~$\mathfrak f$. 
\item[{\it (iii)}] {The test statistic $\tenq{{Q}}^{(n)\pms}_{J_{\mathfrak f}}$ is asymptotically equivalent, under the null $H\n_{0;\text{\,\rm ell}} ({\boldsymbol\beta}^0)$ obtained from $H\n_0 ({\boldsymbol\beta}^0)$ by restricting to elliptical noise and any alternative $H\n_1({\boldsymbol\beta}^0, {\bf B},f)$ where $f$ is elliptical, 
to the test statistic $\utilde{Q}^{(n)\,\text{\rm ell}}_{{J_{\mathfrak f}}}$ based on elliptical ranks and signs. }
\end{compactenum}
\end{Corollary}

 Formal ARE results straightforwardly follow as ratios of standardized shifts or noncentrality parameters; rather than overloading the paper with cumbersome  formulas  we omit explicit expressions. Chernoff-Savage inequalities similar to those obtained in the two-sample case by \cite{deb2021} also follow for the tests based on normal or van der Waerden scores (see Section~\ref{standscore}); in view of Corollary~5.2 {\it (iii)},  these inequalities, which are  limited to the family of elliptical densities, coincide with those  in \cite{HallinPain2002a} and \cite{HallinPain2005}.

  \subsection{Two particular cases}\label{sec5.3}

 In this section, we provide explicit forms of the test statistic for the two-sample and MANOVA problems. Because of their simplicity  
  and practical value (see Section~\ref{standscore}), we concentrate on the case~\eqref{Qpmell} of spherical scores, from which the general case \eqref{Qpm} is easily deduced---essentially, by substituting ${\bf J}\Big(\frac{R\n_{i;{\pms}}}{n_R+1}{{\bf S}\n_{i;{\pms}}}\Big) $ for~$J\Big(\frac{R\n_{i;{\pms}}}{n_R+1}\Big) {{\bf S}\n_{i;{\pms}}}$.

\subsubsection{Center-outward rank tests for  two-sample location}\label{sec-2slt}
An important particular case  is the two-sample location model,   where~$n=n_1+n_2$ and~\eqref{MRegr} holds  
with covariates of the form 
${\mathbf C}\n = (
{\bf 1}_{n_1}^{\prime}, {\bf 0}_{n_2}^{\prime})^{\prime}$  (with~${\bf 1}_{n_1}$   an $n_1$-dimen\-sional column vector of ones, 
$ {\bf 0}_{n_2}$  an~$n_2$-dimensional column vector of zeros); the para\-meter~$\boldsymbol{\beta}=(\beta_{11},\dots,\beta_{1d})^{\prime}$ here is a $d$-dimensional row vector. 
 The objective is to test the null hypothesis $H_0: \ \boldsymbol{\beta}=\mathbf{0}_d$ under which the distributions of~${\bf Y}\n_1,\ldots,{\bf Y}\n_{n_1}$ and 
${\bf Y}\n_{n_1+1},\ldots,{\bf Y}\n_{n}$ coincide. Elementary computation yields  
 $\bar{c}\n=n_1/n,\ \  {V}_{\bf c}\n=n_1n_2/n^2, \ \ \text{ and}~K_{\bf c}=1.$ 
 If the regular  grid  $\mathfrak{G}_n$ is chosen such that $\Vert \sum_{s=1}^{n_S}{\bf s}^{n_S}_s\Vert ={\bf 0}$ (which is always possible in view of  Section~\ref{RSsec}),~$\sum_{i=1}^n J\Big(\frac{R\n_{i;{\pms}}}{n_R+1}\Big){\bf S}\n_{i;{\pms}} ={\bf 0}$ and the test statistic \eqref{Qpmell}  takes the simple form\vspace{-1mm} 
\begin{equation}\label{2sample}
\tenq{Q}^{(n)\pms}_J
=
\frac{nd}{n_1n_2\!\int_0^1 \! J^2(u){\rm d}u}
\left\|
\sum_{i=1}^{n_1}  J\Big(\frac{R\n_{i;{\pms}}}{n_R+1}\Big) {{\bf S}\n_{i;{\pms}}} \right\|^2;
\end{equation}
else, a centering term $\frac{n_1}{n}\sum_{i=1}^n J\Big(\frac{R\n_{i;{\pms}}}{n_R+1}\Big){\bf S}\n_{i;{\pms}} $ is to be subtracted. 
Assumption~\ref{Ass31} {\it (iii)} requi\-res~$\lim_{n\to\infty} n \, {\min\{n_1,n_2\}}/{\max\{n_1,n_2\}} = \infty,  
$
which holds whenever both $n_1$ and $n_2$ tend to infinity.  
 Under this condition and Assumptions~\ref{Ass21}, 
 with   ${\rm P}\in{\cal P}_d$, $\tenq{Q}^{(n)\pms}_J$ is,  under $H_0$, asymptotically $\chi^2$   with~$d$ degrees of freedom and the null hypothesis  can be rejected at asymptotic level $\alpha$ whenever~$\tenq{Q}^{(n)\pms}_J$   exceeds the $(1-\alpha)$ quantile of a $\chi^2_d$ distribution.  
 
{ Noncentrality parameters, in this special case, are particularly simple. Consider the sequence of  alternatives under which the error density is $f$ and $ \boldsymbol{\beta}_n =   n^{-1/2}\mathbf{s}+ o(n^{-1/2})$ with  $\mathbf{s} \in \mathbb{R}^d$. Assume that $n_1/n \to
p\in (0,1)$ as $n\to\infty$. 
   The joint asymptotic normality of vec$\sum_{i=1}^{n_1}  J\Big(\frac{R\n_{i;{\pms}}}{n_R+1}\Big) {{\bf S}\n_{i;{\pms}}}$ and the log-likelihood ratio   follows from a routine application of the Wold-Cram\' er device;  Le Cam's third lemma then readily provides {the asymptotic noncentrality parameter $q = \left( d / \int_0^1 J^2(u) du \right)\Vert\boldsymbol{\mu} \Vert^2$} where, denoting by $\mathbf{Q}_{\pm}$ the center-outward quantile function of the error, 
$
  \boldsymbol{\mu}=
     \left(p(1-p)\right)^{1/2}
   \int_{\mathbb
    S_d} \mathbf{J}(\tu){\boldsymbol\varphi}^\prime_f(\mathbf{Q}_{\pm}(\tu))\,
  \mathrm{d}{\rm U}_d(\tu)\, \mathbf{s}.
$}

\subsubsection{Center-outward rank tests for MANOVA}
Another important special case of   model~\eqref{MRegr}  is the multivariate $K$-sample location or MANOVA model.  The observation here decomposes into $K$ samples, with respective sizes~$n_1,\dots,n_{K}$ and~$n=\sum_{k=1}^K n_k$. Precisely,  
 ${\bf Y}\n=:\big({\bf Y}^{(n;1)\top},\ldots,{\bf Y}^{(n;k)\top},\ldots,{\bf Y}^{(n;K)\top}\big)^\top$ 
with
\[{\bf Y}^{(n;k)}= 
\left(\begin{array}{cccc}
 Y\n_{k;11}& Y\n_{k;12}&\ldots&Y\n_{k;1d}\\ 
\vdots &\vdots&&\vdots \\ 
Y\n_{k;n_k1}& Y\n_{k;n_k2}&\ldots&Y\n_{k;n_kd}
\end{array}\right)
\]
and  \eqref{MRegr}  holds with the   matrix of covariates 
\[
{\mathbf C}\n= \begin{pmatrix} {\bf 1}_{n_1}& {\bf 0}_{n_1} & \dots & {\bf 0}_{n_1}\\
{\bf 0}_{n_2}& {\bf 1}_{n_2} & \dots & {\bf 0}_{n_2}\\
\vdots & \vdots & \cdots &\vdots\\
{\bf 0}_{n_K}& {\bf 0}_{n_K} & \dots & {\bf 0}_{n_K}\end{pmatrix}=
\begin{pmatrix} 
\td\n_{11} & \td\n_{12} & \dots & \td\n_{1,K-1}\\
\td\n_{21}&\td\n_{22}&\dots &\td\n_{2,K-1}\\
\vdots & \vdots & \cdots &\vdots\\
\td\n_{K1}&\td\n_{K2}&\dots&\td\n_{K,K-1},
\end{pmatrix} 
\]
where 
 $\td\n_{ij}\!\!={\bf 1}_{n_i} I[i=j]$, 
 $i=1,\dots,K$ and $j=1,\dots,{K\!-1}$.  The null hypothesis   
  is the hypothesis of no treatment effect $H_0: \ \boldsymbol{\beta}={\mathbf 0}_{(K-1)\times d}$. 
 
Letting ${\bf v}\n:=(n_1/n,\dots,n_{K-1}/n)\pr$,   the matrix ${\bf V_c}\n $  in Assumption \ref{Ass31} takes the form~$
 {\bf V_c}\n 
= \mathrm{diag}\{ {\bf v}\n\} -  { {\bf v}\n}{ {\bf v}\n}\pr,
$ 
 where $\mathrm{diag}\{ {\bf v}\n\}$ stands for the diagonal matrix with  diagonal entries $ {\bf v}\n$. 
If the regular  grid  $\mathfrak{S}_n$ is chosen such that~$\Vert \sum_{s=1}^{n_S}{\bf s}^{n_S}_s\Vert =0$ and~$
(\mathbf{V}_{\mathbf{c}}^{(n)})^{-1/2}$ is substituted for its limit $\mathbf{K}_{\mathbf{c}}^{(n)}$,   the test statistic~\eqref{Qpmell}  simplifies into\vspace{-2mm} 
\[
{\tenq{Q}}^{(n)\pms}_J = \frac{d}{\int_{0}^1 J^2(u){\rm d}u } \sum_{k=1}^K \frac{1}{n_k}\left\|\sum_{i=n_1+\ldots+n_{k-1}+1}^{n_1+\ldots+n_k} \hspace{-5mm}J\Big(\frac{R\n_{i;{\pms}}}{n_R+1}\Big) {{\bf S}\n_{i;{\pms}}} \right\|^2.
\]

Assumption  \ref{Ass31}{\it (iii)}  is satisfied as soon as 
  $\lim_{n\to\infty}  \min\{n_1,\dots,n_K\}\to\infty.$  
 Assuming moreover that 
$0<\lim\inf_{n\to\infty} n_k/n \leq \lim\sup_{n\to\infty} n_k/n  <1$ 
 for  $1\leq k\leq K$,    the limit matrix~${\mathbf R_c}$ is positive definite\footnote{This limit possibly\vspace{-2mm} can exist along subsequences, with asymptotic statements modified accordingly. For the sake of simplicity, we do not include this in subsequent results. } and  Assumption~\ref{Ass31}{\it (ii)}  is satisfied as well. Then, under the null  hypothesis  of no treatment effect,
    ${\tenq{Q}}^{(n)\pms}_J$ is asymptotically chi-square with~$(K-~\!1)d$ degrees of freedom and the test rejecting $H_0$ whenever ${\tenq{Q}}^{(n)\pms}_J$ exceeds the corresponding~$(1-~\!\alpha)$ quantile has asymptotic level $\alpha$ irrespective of the actual error distribution~${\rm P}\in{\mathcal P}_d$.  This test is a multivariate generalization of the well-known univariate rank test  for  $K$-sample equality of location (the univariate one-way ANOVA hypothesis of  no treatment effect), see \cite[p.170]{HajekSidak}. Note that, for $K=2$,   ${\tenq{Q}}^{(n)\pms}_J $  coincides with the two-sample test statistic obtained in Section~\ref{sec-2slt}.

\section{Choosing a score function} Section~\ref{Rpmsec}  allows us to construct, based on any 
$\bf J$ or $J$  satisfying Assumption~\ref{Ass21} (either with~\eqref{monass} or \eqref{Jsphcond}),  strictly distribution-free center-outward rank tests of the null
 hypothesis $H\n_0({\boldsymbol\beta}^0)$   under    which~${\boldsymbol\beta}=~\!{\boldsymbol\beta}^0$ while the intercept  ${\boldsymbol\beta}_0$ and the error distribution ${\rm P}\in{\mathcal P}_d$ 
remain unspecified. All these  tests, however, depend on a score function to be selected by the practitioner. Some will favor simple scores of the spherical type (see Section~\ref{standscore}); others may want to base their choice on efficiency considerations (see Section~\ref{effscore}). 

\subsection{Standard score functions}\label{standscore}

  Popular choices   are the spherical sign test, Wilcoxon  and van der Waerden scores. Let us describe them, in more details,  in the particular case of the two-sample  problem.

The two-sample sign test is based on   the  degenerate score  $J_{\text{\tiny sign}}(r)  :=1$ for $r \in [0,1)$;  using the fact that $\sum_{i=1}^{n} {\bf S}_{i;{\pms}}^{(n)}={\bf 0}$, 
 one gets for \eqref{2sample}, with the notation of Section~\ref{sec-2slt}, the very simple test statistic
\[
  \tenq{{  Q}}^{(n)\pms}_{\text{\tiny sign}} = \frac{nd}{n_1n_2}\left\|\sum_{i=1}^{n_1} {\bf S}_{i;{\pms}}^{(n)} \right\|^2.
\]

The choice $J_{\text{\tiny Wilcoxon}}(r):=r$ similarly characterizes the Wilcoxon two-sample test: noting that~$\sum_{i=1}^{n} R\n_{i;{\pms}}{\bf S}_{i;{\pms}}^{(n)}\! ={\bf 0}$ holds if $\sum_{i=1}^{n} {\bf S}_{i;{\pms}}^{(n)}\! ={\bf 0}$ and that $\int_0^1 r^2 {\rm d}u\! =~\!\!1/3$, this yields\vspace{-1mm}
\[
 \tenq{Q}^{(n)\pms}_{\text{\tiny Wilcoxon}} = \frac{3nd}{n_1n_2{(n_R+1)^2}}\left\|\sum_{i=1}^{n_1} R\n_{i;{\pms}}{\bf S}_{i;{\pms}}^{(n)} \right\|^2.
\]
For $d=1$ the test based on $ \tenq{{Q}}^{(n)\pms}_{\text{\tiny Wilcoxon}}$ is asymptotically equivalent to the classical univariate  two-sided two-sample Wilcoxon test.
 
As for the  two-sample van der Waerden test, it is based on the Gaussian or van der Waerden scores $J_{\text{\tiny vdW}}(r) := \big(\Psi_d^{-1}(r)\big)^{1/2}$, where $\Psi_d$ denotes  the cumulative distribution function of   a chi-square variable with $d$ degrees of freedom.
 Clearly $\int_0^1 J^2_{\text{\tiny vdW}}(r) {\rm d}r = \int_0^\infty x {\rm d} \Psi_d(x) = d$    and, provided that~$\sum_{i=1}^{n} {\bf S}_{i;{\pms}}^{(n)}=0$, 
$ \sum_{i=1}^{n}\Big( \Psi_d^{-1}\Big(\frac{R\n_{i;{\pms}}}{n_R+1}\Big)\Big)^{1/2} {\bf S}_{i;{\pms}}^{(n)}=0$.  Hence, the van der Waerden center-outward rank test statistics takes the~form 
\[
  \tenq{{Q}}^{(n)\pms}_{{\text{\tiny vdW}}} = \frac{n}{n_1n_2}\left\|\sum_{i=1}^{n_1} \Big( \Psi_d^{-1}\Big(\frac{R\n_{i;{\pms}}}{n_R+1}\Big)\Big)^{1/2} {\bf S}_{i;{\pms}}^{(n)} \right\|^2.
\]

\subsection{Score functions and efficiency}\label{effscore} 

The  {test} statistics in Section~\ref{standscore} offer the advantage of a structure  paralleling the structure of the numerator of the  classical Gaussian $F$ test---basically substituting, in the latter,~$ {\bf S}_{i;{\pms}}^{(n)}\vspace{-1mm}$ (sign test scores),  $R\n_{i;{\pms}} {\bf S}_{i;{\pms}}^{(n)}\vspace{-0mm}$ (Wilcoxon scores), or~$\Big( \Psi_d^{-1}\Big(\frac{R\n_{i;{\pms}}}{n_R+1}\Big)\Big)^{1/2}{\bf S}_{i;{\pms}}^{(n)}$ (van der Waerden scores) for the sphericized residuals \eqref{Zell} (the computation of which,  moreover, requires  the specification of ${\boldsymbol\beta}_0$ or its consistent estimation, something   center-outward ranks and signs do not need in view of their shift-invariance) and adopting the adequate standardization.

The choice of a score function  also can be  guided by efficiency considerations, selecting~$\bf J$ in relation to some reference  distribution under which efficiency is to be attained.  This, in the univariate case, yields the normal (van der Waerden),  Wilcoxon or sign test scores, achieving efficiency under  Gaussian, logistic, or double exponential reference densities; as we shall see,~$\tenq{{\bf Q}}^{(n)\pms}_{\text{\tiny sign}} $ and $  \tenq{{\bf Q}}^{(n)\pms}_{{\text{\tiny vdW}}}$ similarly achieve efficiency at spherical exponential and Gaussian reference distributions. Due to the fact that the density $f^\star_{d;{\mathfrak f}}$ of the modulus of a spherical logistic fails to be logistic for $d>1$,  the Wilcoxon test based on  $\tenq{{\bf Q}}^{(n)\pms}_{\text{\tiny Wilcoxon}} \vspace{1mm}$, however, does not enjoy efficiency under spherical logistic; this is also the case of the   elliptical rank tests based on Wilcoxon scores in \cite{HallinPain2002b, HallinPain2002a, HallinPain2005}.

In the same spirit, one could  contemplate the idea of achieving, based on  center-outward rank tests,  efficiency at some selected reference distribution~$ 
{\rm P}^{\boldsymbol\varepsilon}_0$ in~${\cal P}_d$
  (with density $f_0^{\boldsymbol\varepsilon}$ and center-outward distribution function~${\bf F}_{0;\pms}^{{\boldsymbol\varepsilon}}$ satisfying the adequate regularity assumptions).  
 Indeed, it follows from Proposition~\ref{asNpm} that efficiency under~${\rm P}^{{\boldsymbol\varepsilon}}_0$ 
 can be achieved by  a test based  on the test statistic $\tenq{{\bf Q}}^{(n)\pms}_{\bf J}$ given in~\eqref{Qpm} with score 
 ${\bf J}={\boldsymbol \varphi}_{f^{\boldsymbol\varepsilon}_0}\circ\big({\bf F}^{\boldsymbol\varepsilon}_{0;\pms}\big)^{-1}$. 
 This, however, raises two problems. First, in order for  ${\boldsymbol\varphi}_{f^{\boldsymbol\varepsilon}_0}$ to be  analytically computable, the distribution~${\rm P}_0^{\boldsymbol\varepsilon}$ has to be fully specified (up to location and a global scaling parameter), with  closed-form density function~$f^{\boldsymbol\varepsilon}_0$. Second, the corresponding score function~${\bf J}={\boldsymbol \varphi}_{f^{\boldsymbol\varepsilon}_0}\circ\big({\bf F}^{\boldsymbol\varepsilon}_{0;\pms}\big)^{-1}$ also involves the  center-outward quantile function $({\bf F}^{\boldsymbol\varepsilon}_{0;\pms})^{-1}$ for which, except for   a few particular cases (spherical distributions),  no   explicit form is available in the literature.  Once~${\rm P}^{\boldsymbol\varepsilon}_0$ is fully specified, in principle, it can be simulated, and an arbitrarily precise   numerical evaluation of $({\bf F}^{\boldsymbol\varepsilon}_{0;\pms})^{-1}$ can be obtained, to be plugged into ${\bf J}$. This may~be~computa\-tionally heavy, but increasingly efficient algorithms are available in the    domain of  numerical  measure transportation: see, e.g.,  \cite{merigot2011multiscale} or \cite{peyre+c:2019}. 

Now,  choosing a fully specified reference   ${\rm P}^{\boldsymbol\varepsilon}_0$ may be embarrassing---this means, for instance, a skew-$t$ distribution with specified degrees of freedom, shape matrix, and skewness parameter (without loss of generality, location can be taken as~$\bf 0$), a multinormal or elliptical distribution with specified radial density and specified  (up to a positive global factor)  covariance (again, the mean can be taken as~$\bf 0$), 
  {or any other multivariate distribution with fully specified parameters.}
 Fortunately, a full specification of ${\rm P}^{\boldsymbol\varepsilon}_0$ can be relaxed to the specification of a parametric family with parameter~$\boldsymbol\vartheta$, say, such as the family~${\cal P}_{\text{\tiny\rm skew\! {\it t}}}$ of all skew-$t$ distributions with location $\bf 0$ (parameters: a shape matrix and a~$d$-tuple of skewness parameters) or the family~${\cal P}_{\mathfrak f}^{\text{\tiny{\rm ell}}}$ of all elliptical distributions~\eqref{ellf} with radial density~$\mathfrak f$ (parameter: a scatter matrix). The unspecified parameter~$\boldsymbol\vartheta$ of ${\rm P}^{\boldsymbol\varepsilon}_0$ indeed can be replaced, in the numerical evaluation of ${\bf F}^{\boldsymbol\varepsilon}_{0;\pms}\vspace{-1mm}$, with consistent estimated values provided that the estimator $\hat{\boldsymbol\vartheta}$ is measurable with respect to the {\it order statistic}\footnote{The order\vspace{-.5mm} statistic of the $n$-tuple ${\bf Z}_1,\ldots,{\bf Z}_n$ of $d$-dimensional ($d>1$) random vectors  can be defined as any reorde\-ring~${\bf Z}_{(1)},\ldots,{\bf Z}_{(n)}$ generating the $\sigma$-field of permutation-invariant Borel sets of~$\sigma\big({\bf Z}_1,\ldots,{\bf Z}_n\big)$;\vspace{-.5mm} for instance, the one resulting from ordering the observations ${\bf Z}_i$ from smallest to largest first component.} 
  of the residuals~${\bf Z}\n_i\!$. Plugging these estimators into the score~$\bf J$---this {includes} the standardization factor and the numerical evaluation of ${\bf F}^{\boldsymbol\varepsilon}_{0;\pms}$---yields data-driven (order-statistic-driven) scores~${\bf J}\n$; similar   data-driven scores have been proposed in the univariate case by \cite{DJurec00}. Conditionally on the order statistic,  the corresponding test statistic is still distribution-free and its (conditional) critical values yield  unconditionally  correct size. However, these critical values involve the order statistic: the resulting tests therefore no longer are ranks tests but permutation tests.\footnote{A {\it permutation test }is a test\vspace{-1mm} enjoying {\it Neyman $\alpha$-structure} with respect to the sufficient and complete order statistic.}  
 The theoretical properties,  feasibility, and finite-sample performance of this data-driven approach should be explored and   numerically assessed---this is, however, beyond the scope of this paper and we leave it for future research.

In view of this,  no obvious non-spherical  convenient  candidate emerges   as a  reference density  in dimension~$d>1$. The center-outward 
test statistic achieving optimality at the spherical distributions with radial density~$\vspace{-1mm}{\mathfrak f}$ is~$\tenq{{\bf Q}}^{(n)\pms}_{J_{\mathfrak f}}$ with~$J_{\mathfrak f}$ as  in part {\it (iii)} of Corollary~\ref{scorell}. 

\color{black}
\subsection{Affine invariance and sphericization}\label{invSec}

Affine invariance (testing) or equivariance (estimation), in ``classical multivariate analysis,'' is often considered an essential and inescapable property. Closer examination, however, reveals that this particular role of affine transformations is intimately related to the affine invariance of Gaussian and elliptical families of distributions. When Gaussian or elliptical assumptions are relaxed,  affine transformations are losing this privileged role and the relevance of affine invariance/equivariance properties is much less obvious. We refer to Appendix~A.6 for a more detailed discussion of that invariance issue.

\setcounter{equation}{0}

\section{Some numerical results} 

A Monte Carlo simulation study is conducted (Sections~\ref{Sec7.1}--\ref{Sec7.2}) in order to explore the finite-sample performance of our tests. Results  are presented for two-sample location and MANOVA models,  and limited to the  Wilcoxon score function $J(r)=r$;  other choices for~$J$ lead to very similar figures, which we therefore do not report. {The analysis was conducted in R program \cite{R}. The center-outward ranks and signs were computed using the optimal transportation via the so-called Hungarian algorithm implemented in package \texttt{clue}, \cite{clue}. The Supplementary Material contains more details on the exact implementation.}

\subsection{Two-sample location, $d=2$ }\label{Sec7.1}

Consider first the two-sample location problem in dimension $d=2$. Two independent random samples of size $n_1=n_2=n/2$ were generated and 
the two test statistics $\tenq{Q}^{(n)\,\text{\rm{ell}}}_{\text{\tiny Wilcoxon}}$ and~$\tenq{Q}^{(n)\pms}_{\text{\tiny Wilcoxon}} \vspace{-0.5mm}$ (see Sections~\ref{Rpmsec} and~\ref{sec-2slt}) 
were computed. 
 The sample covariance matrix $\widehat{\tSigma}$ was used for the computation of $\tenq{{\boldsymbol\Lambda}}^{(n)\,\text{\rm{ell}}}_{\text{\tiny Wilcoxon}}$ and the elliptical or Mahalanobis ranks and signs.  
 
  Rejection frequencies were computed for the following error densities: 
\begin{compactenum}
\item[{\it (a)}] a centered bivariate normal distribution  with  unit variances and  correlation $\rho=1/4$;
\item[{\it (b)}] a centered bivariate $t$-distribution  with   the same scaling matrix as in {\it (a)} and $\nu$  degree of freedom,  $\nu=1$ (Cauchy) and $\nu=3$;\footnote{The bivariate\vspace{-.5mm} $t$-distribution with $m$ degrees of freedom and scaling matrix ${\bf A}^\prime{\bf A}$ is the one  defined in Example 2.5 of  \cite{fang17} as the distribution of a random vector\vspace{-.5mm} ${\boldsymbol\xi}:= {\boldsymbol\mu} + {\bf A}^\prime {\boldsymbol   \zeta}\sqrt{m}/\sqrt{s}$\linebreak  where~${\boldsymbol \zeta}\sim\mathcal{N}_2({\bf 0},{\bf I}_2)$\vspace{-.5mm}   and $s \sim\chi ^2_m$, independent of ${\boldsymbol \zeta}$---not to be confused with the elliptical  distribution with Student radial density $\mathfrak f$.} 
\item[{\it (c)}] a mixture, with weights
$w_1=1/4$ and $w_2=3/4$,  of two bivariate normal distributions with means $\tmu _1 =(3/4,0)\pr$ and $\tmu_2 =(-1/4,0)\pr$ and  covariance matrices 
\[
\tSigma_1= \begin{pmatrix} 1&\; 2/3\\\; 2/3 & 1\end{pmatrix}\quad\text{ and }\quad \tSigma_2= \begin{pmatrix} 1&-2/3\\ -2/3 & 1\end{pmatrix}, \vspace{-4mm}
\]
 respectively;
\item[{\it (d)}] a mixture, with weights $w_1=1/4$ and $w_2=3/4$,  of two bivariate $t_1$ (Cauchy) distributions centered at $\tmu _1 =(3/4,0)\pr$ and $\tmu_2 =(-1/4,0)\pr$,   with the same scaling matrices $\tSigma_1$ and $\tSigma_2$ as in {\it (c)};
 \item[{\it (e)}] a  ``U-shaped'' mixture, with weights $w_1=1/2$, $w_2=1/4$, and $w_3=1/4$,  of three bivariate normal distributions, $\mathcal{N}_2(\tmu_1,\tSigma_1)$, $\mathcal{N}_2(\tmu_2, \tSigma_2)$, and $\mathcal{N}_2(\tmu_3,\tSigma_3)$ where \vspace{-2mm}
\[
\tmu_1=(0,0)\pr,\quad \tmu_2=(-3,1)\pr, \quad \tmu_3=(3,1)^\prime ,\vspace{-2mm}
\]
and\vspace{-2mm}
\[
\tSigma_1= \begin{pmatrix} 2\; &0\\ 0\; &  {1}/{8}\end{pmatrix}, \quad \tSigma_2=\begin{pmatrix} {1}/{2}&- {1}/{3}\\  - {1}/{3} & {1}/{2}\end{pmatrix}, \quad \tSigma_3=\begin{pmatrix}  {1}/{2}\;& {1}/{3}\\  1/3\; & 1/2\end{pmatrix};
\] 
\item[{\it (f)}] an ``S-shaped'' mixture, with equal weights $w=1/3$,   of three bivariate normal distributions, $\mathcal{N}_2(\tmu_4,\tSigma_4)$, $\mathcal{N}_2(\tmu_5, \tSigma_5)$, and $\mathcal{N}_2(\tmu_6,\tSigma_4)$  where\vspace{-2mm}
\[
\tmu_4=(-9/2,-1/2)^\prime,\quad\tmu_5=(0,-1/2)^\prime,\quad\tmu_6=(9/2,1)^\prime,
\]
and\vspace{-2mm}
\[
\tSigma_4= \begin{pmatrix} 3/2&-\sqrt{3/8}\\-\sqrt{3/8}&1 \end{pmatrix}, \quad \tSigma_5= \begin{pmatrix} 3/2&\sqrt{3/8}\\ \sqrt{3/8}&1 \end{pmatrix}, \quad \tSigma_6= \begin{pmatrix} 3/2&-\sqrt{3/8}\\-\sqrt{3/8}&1 \end{pmatrix};
\]
\item[{\it (g)}] a skew-$t$-distribution with $\nu$ degrees of freedom, $\nu=1$ and $3$, with  skewness para\-meter~$\boldsymbol{\alpha}=(5,-3)^\prime$,  scaling matrix $\tSigma_7= \begin{pmatrix} 1\; &-1/2\\ -1/2\; & 1\end{pmatrix}$, 
and   location $\boldsymbol{\xi}={\bf 0}$. 
\end{compactenum}

Mixture error densities   naturally appear in the context of hidden heterogeneities due, for instance, to omitted covariates; as for asymmetries, they are likely to be the rule rather than the exception.  Samples of size $200$ from the Gaussian mixtures {\it (c)}, {\it (e)}, and~{\it (f)} and  the skew-$t$ distribution with~$3$ degrees of freedom {\it (g)} are shown in Appendix~A.7.1,~Figure~A.3.

To investigate finite-sample performance, a first sample was generated from one of the  distributions {\it (a)--(g)},   a second one   from the same distribution shifted by the vector $(\delta,\delta)^{\prime}$ for $\delta\in[0.00,0.24]$. Three  sample sizes~$n_1=n_2=50,\, 200$, and~450  (hence, $n=100,\, 400$, and  900) were considered, yielding three groups of curves (from light gray to black, colors in the online version).   The regular grids   $\mathfrak{G}_n$  for computation of the center-outward ranks and signs are constructed with~$n_S=n_R=10$ for~$n=100$,  $n_S=n_R=20$ for $n=400$, and  $n_S=n_R=30$  for $n=900$.  
Each simulation was replicated  $N=1000$ times and the empirical size and power of the test were computed for~$\alpha=0.05$. The resulting rejection frequencies show the dependence of the power on~$\delta$;  they
 are provided in Figures~\ref{fig1}--\ref{fig4}.  For the sake of comparison, we also provide the power of Hotelling's classical two-sample   test. \vspace{-6mm}. 
 \begin{figure}[h!]
\begin{tabular}{lll}
\hspace{-5mm}\includegraphics[width=0.36\textwidth, height=0.44\textwidth]{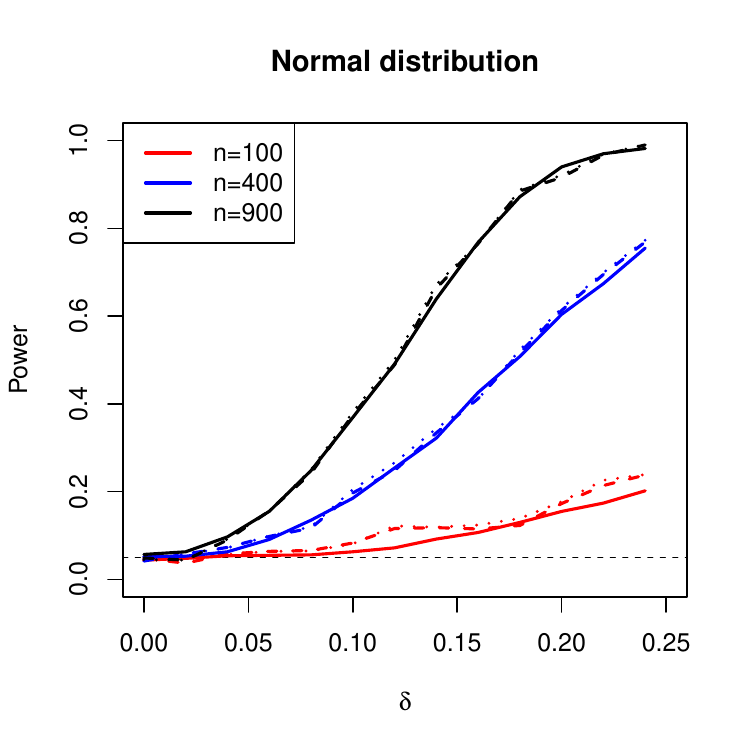} &\hspace{-8mm}
\includegraphics[width=0.36\textwidth, height=0.44\textwidth]{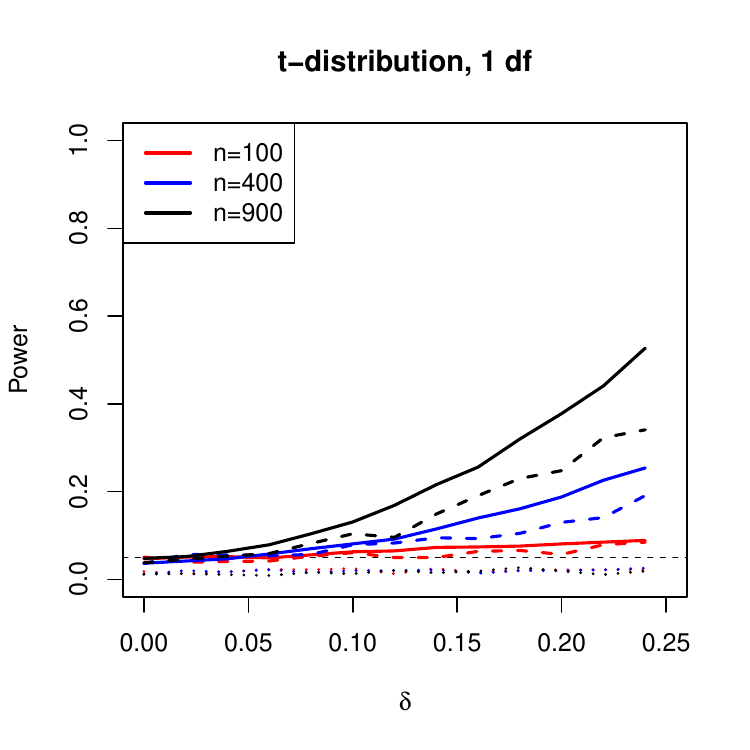} &\hspace{-8mm}
\includegraphics[width=0.36\textwidth, height=0.44\textwidth]{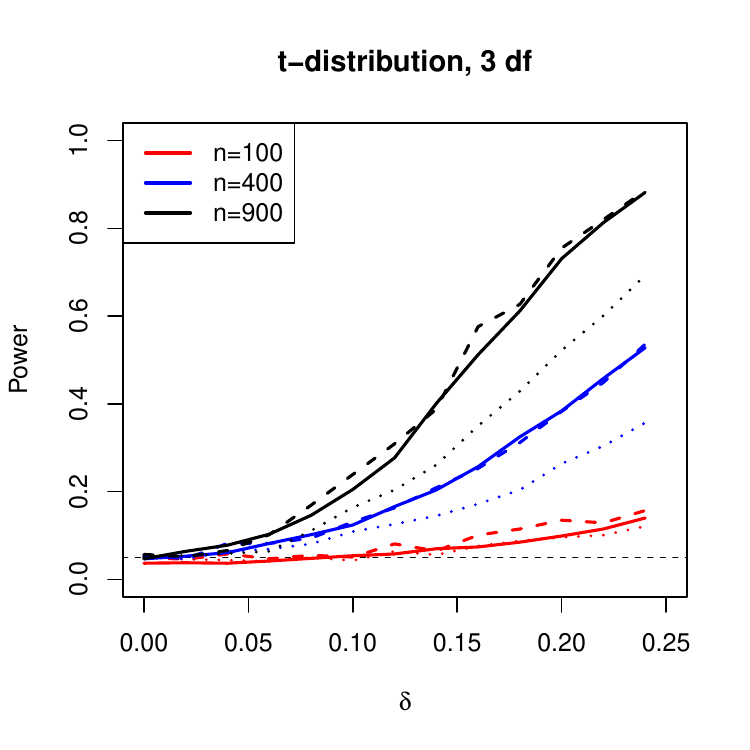} \vspace{-8mm}
\end{tabular}
 \caption{\small Empirical powers of two-sample location tests based on\vspace{-.5mm}  
 the Wilcoxon center-outward rank statistic (solid line),
 the Wilcoxon 
 elliptical\vspace{-.5mm} rank  statistic  
(dashed line), and  Hotelling's  two-sample test (dotted line), as  functions of the shift $\delta$ under\vspace{-.5mm} bivariate normal  and elliptical Student (1 and~3 degrees of freedom) error densities; sample sizes~$n_1=n_2=50$ (red), $200$ (blue), and  $450$ (black).\vspace{-4mm}}\label{fig1}
\end{figure}

Figure~\ref{fig1} displays the empirical power curves for the elliptical distributions {\it (a)} and~{\it (b)}. The results 
  for the normal distribution are very similar for the three tests:   rank-based tests (Wilcoxon scores), thus, are no less powerful than the optimal Hotelling test.  
  As expected,   Hotelling  crashes under  the $t_1$ distribution, while  
the Wilcoxon elliptical  test, although based on the sample covariance matrix, performs surprisingly well  (the robustness  of ranks offsets infinite variance).  
The   tests based  on $\tenq{Q}^{(n)\,\text{\rm{ell}}}_{\text{\tiny Wilcoxon}}$ and $\tenq{Q}^{(n)\pms}_{\text{\tiny Wilcoxon}}\vspace{0.2mm}$ both outperform   Hotelling   also for the $t$-distribution with 3 degrees of freedom. The conclusion is that center-outward rank tests perform equally well as   elliptical  rank tests under elliptical~densities.\vspace{-0mm} 

\begin{figure}[h!]
\centering
\includegraphics[width=0.36\textwidth, height=0.44\textwidth]{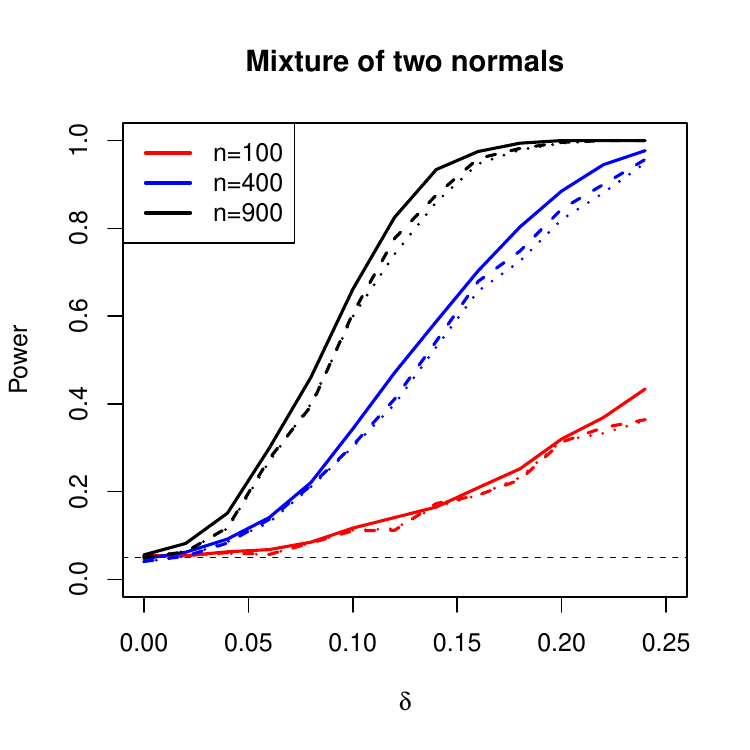}  
\includegraphics[width=0.36\textwidth, height=0.44\textwidth]{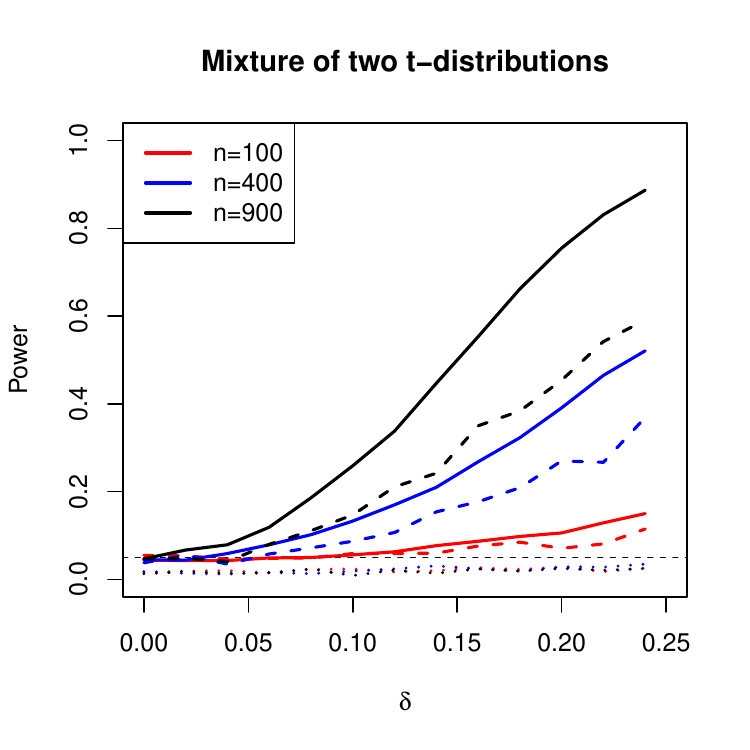} \vspace{-8mm}
\caption{\small Empirical powers of two-sample location tests based on\vspace{-.5mm}   
 the  Wilcoxon center-outward rank statistic (solid line),
 the 
Wilcoxon elliptical\vspace{-.5mm} rank test statistic  
(dashed line), and Hotelling's  two-sample  test (dotted line), as  functions\vspace{-.5mm} of the shift $\delta$, for   the mixtures of two normal   (left panel) and two $t_1$ error densities (right panel), respectively; sample\vspace{-.5mm} sizes~$n_1=n_2=50$ (red), $200$ (blue), and  $450$ (black).\vspace{-3mm}}\label{fig2}
\end{figure}

\begin{figure}[t!]
\centering
\begin{tabular}{cc}
\includegraphics[width=0.36\textwidth, height=0.44\textwidth]{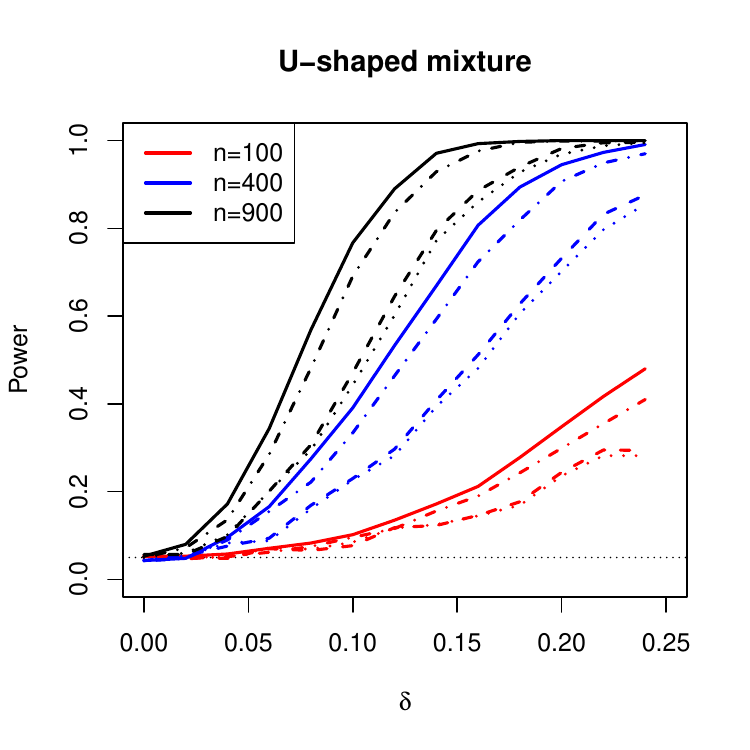}  &
\includegraphics[width=0.36\textwidth, height=0.44\textwidth]{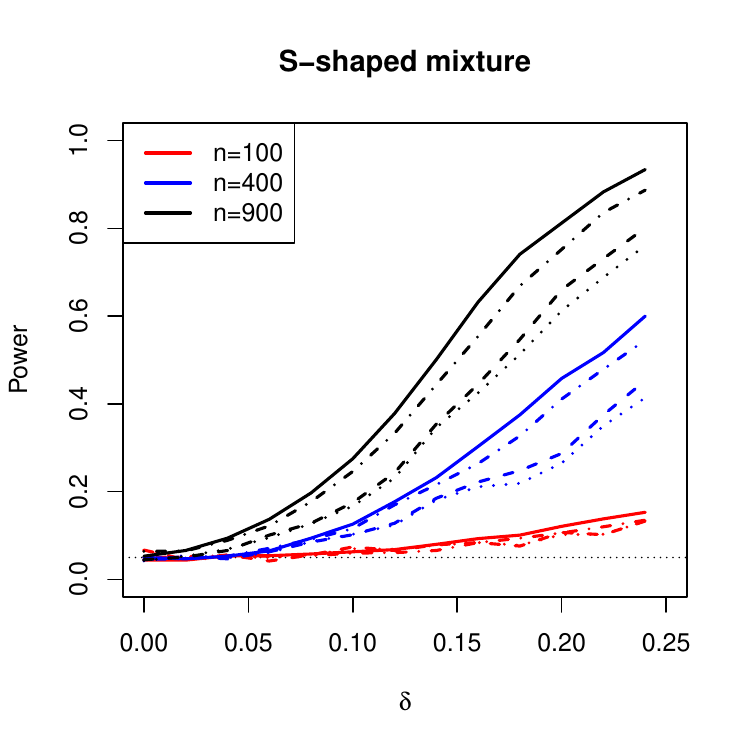} \vspace{-7mm} \\ 
\includegraphics[width=0.36\textwidth, height=0.44\textwidth]{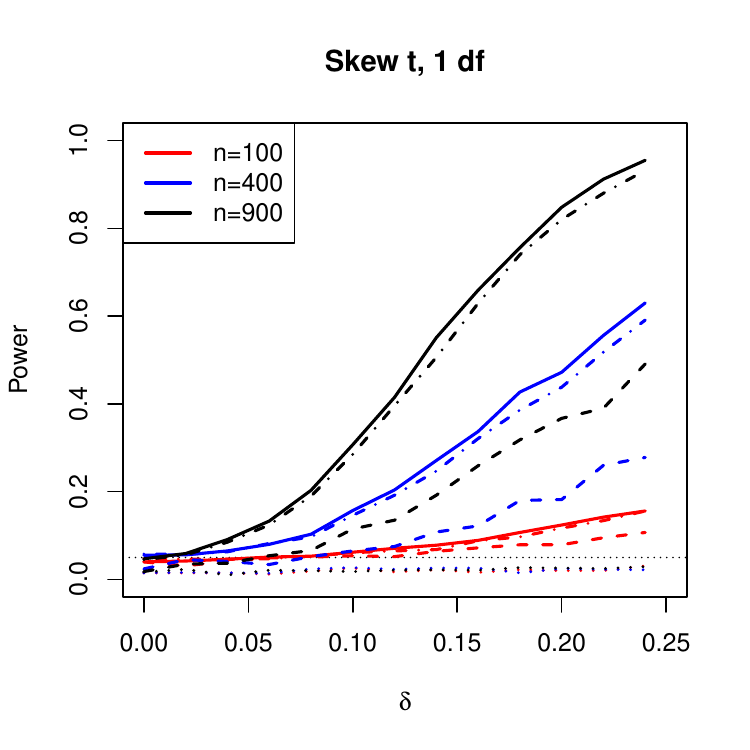} &
\includegraphics[width=0.36\textwidth, height=0.44\textwidth]{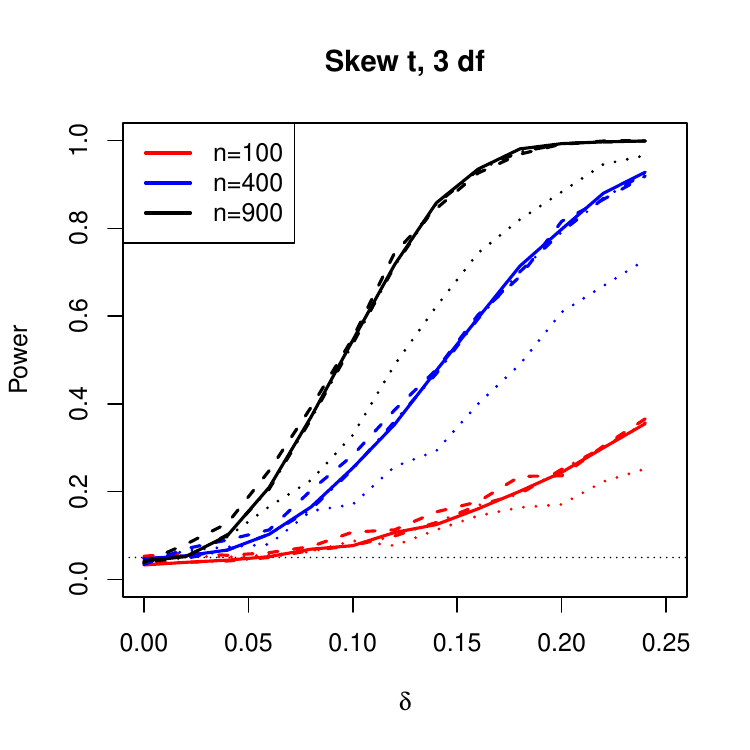} \vspace{-8mm}
\end{tabular}
\caption{\small Empirical powers of two-sample location tests based on \vspace{-.5mm}  
 the Wilcoxon center-outward rank statistic (solid line), the Wilcoxon center-outward rank \vspace{-.5mm}    statistic computed from linearly sphericized  residuals (dot-dashed line), 
  the Wilcoxon   elliptical rank test statistic\vspace{-.5mm}  
(dashed line),  and Hotelling's two-sample test (dotted line),  
as  functions of the shift $\delta$   for the "U-shaped"  \vspace{-.5mm}  (upper left panel) and the "S-shaped"  (upper right panel) mixtures\vspace{-.5mm}  of three normal error densities,  and 
 skew-$t$ error densities with\vspace{-.5mm}  $\nu=1.1$   (bottom left panel) and $\nu=3$  (bottom right panel) degrees of freedom, respectively; sample sizes~$n_1=n_2=50$ (red), $200$ (blue), and  $450$ (black). \vspace{-3mm}}\label{fig4}
\end{figure}

The remaining distributions {\it (c)--(g)} are non-elliptical ones. Results for the mixtures {\it (c)} and {\it (d)} are shown in Figure~\ref{fig2}. For the mixture {\it (c)} of two normals, the results obtained for the three   tests are still quite similar, but the center-outward rank  test based on~$\tenq{Q}^{(n)\pms}_{\text{\tiny Wilcoxon}} $, in general,  yields the largest power. 
For the mixture {\it (d)} of two $t_1$ (Cauchy) distributions, the Hotelling test fails miserably and
the center-outward  rank test very clearly outperforms the elliptical   rank test for all   sample sizes.
 Figure~\ref{fig4} provides the results for the mixtures {\it (e)--(f)} and the skew-$t$-distribution~{\it (g)},  respectively. The power curve for the test statistic~$\tenq{Q}^{(n)\pms}_{\text{\tiny Wilcoxon}} $ computed from the linearly sphericized residuals (using the sample mean and the sample covariance matrix as estimators of location and scatter) is added as a dot-dashed line.
 In all these plots, the 
center-outward rank test statistic leads to the largest power. Note that the linear sphericization of the residuals, which makes the test affine-invariant,  may noticeably deteriorate the power (see the discussion in Section~\ref{invSec} and Appendix~A.6).

\begin{figure}[ht!]
\centering
\includegraphics[width=0.36\textwidth, height=0.44\textwidth]{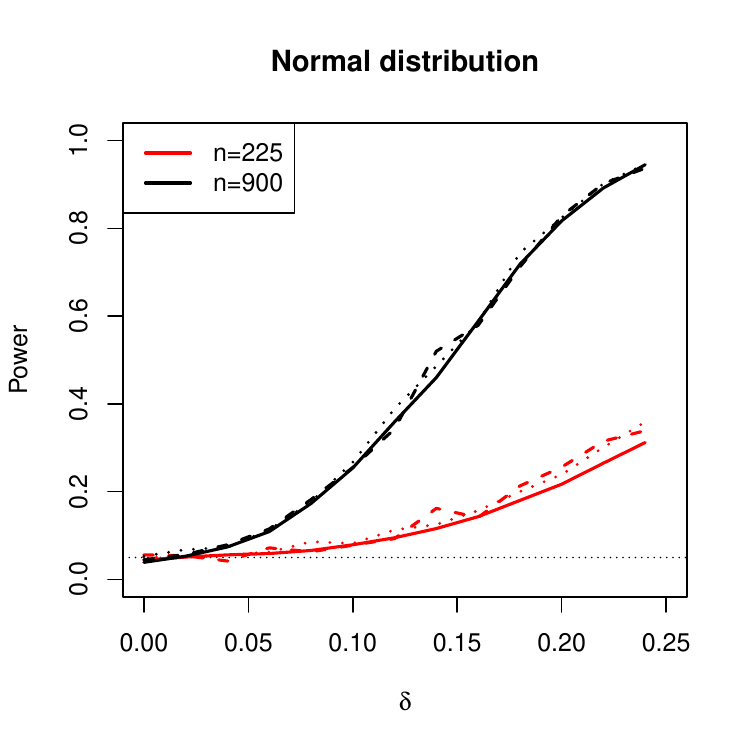} 
\includegraphics[width=0.36\textwidth, height=0.44\textwidth]{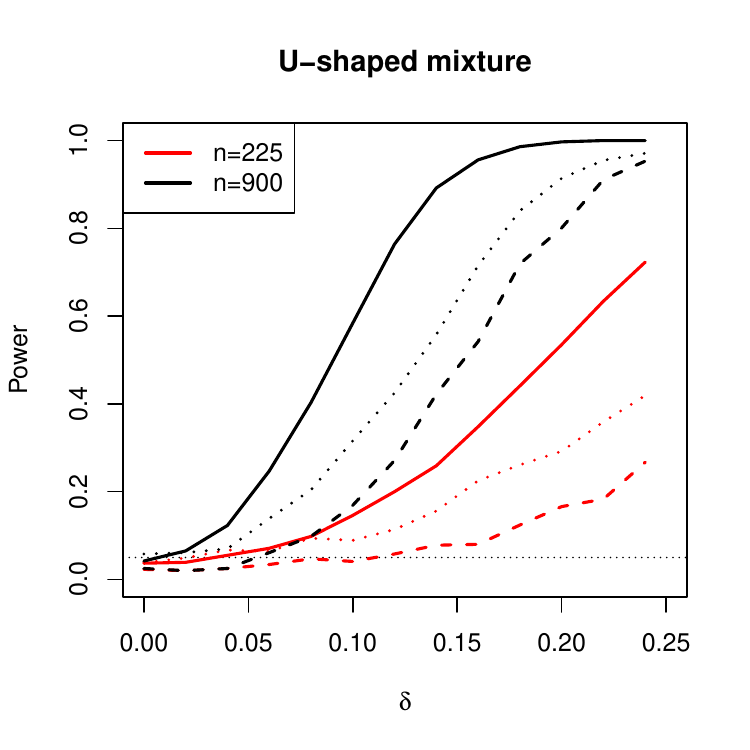}  \vspace{-8mm}
\caption{\small Empirical powers of  MANOVA tests based\vspace{-.5mm} on the Wilcoxon center-outward rank statistic (solid line), the Wilcoxon  elliptical rank test statistic \vspace{-.5mm} 
(dashed line), Pillai's  test (dotted line), and  {Roy's test (dashed-dotted line)} 
 as  functions of the shift $\delta$, for the normal distribution\vspace{-.5mm} (left panel) and the U-shaped mixture of three normals (right panel); the sample sizes   \vspace{-.5mm}  are~$n_1=n_2=n_3=75$ (red)   and~$300$ (black).  \vspace{-3mm}
 }\label{fig-manova}
\end{figure}

\subsection{One-way MANOVA, $d=2$}\label{Sec7.2}
The performance of center-outward rank tests is very briefly studied here  for one-way MANOVA with $K=3$ groups, still for $d=2$. Two random samples were generated from the distribution {\it (a)}  (Gaussian) or {\it (e)} (U-shaped mixture of three Gaussians), as described in Section \ref{Sec7.1}, and the third sample was drawn from the same distribution shifted  by the vector $(\delta,\delta)^{\prime}$  for $\delta\in[0.00,0.24]$.  A balanced design with  groups of size $n_1=n_2=n_3=75$ (hence~$n=225$)  and $n_1=n_2=n_3=300$ (hence~$n=900$)  was considered. For $n=225$, the grid  $\mathfrak{G}_n$     is constructed with $n_R=n_S=15$; for $n=900$, we set~$n_R=n_S=30$.  As in Section~\ref{Sec7.1},  the results are presented for the Wilcoxon scores~$J(r)=r$ only---other choices lead to very similar conclusions.

Rejection frequencies are plotted in Figure~\ref{fig-manova} for the center-outward rank test based on~$\tenq{Q}^{(n)\pms}_{\text{\tiny Wilcoxon}} $  (solid line), the elliptical rank test statistic 
$\tenq{Q}^{(n)\,\text{\rm{ell}}}_{\text{\tiny Wilcoxon}}$ (dashed line), 
 the Pillai trace test based on an approximate F-distribution  (dotted line), {and the Roy test (dashed-dotted line); see the Supplementary Material for implementation.  Under normal density,  all the tests} perform very similarly.  For the non-elliptical mixture distribution,  however, 
 the center-outward rank test achieves  sizeably larger power than all {other ones.}  
 Further simulations yielding, in dimension $d=6$, similar conclusions, are provided in~Appendix~A.7. 

\subsection{An empirical illustration}\label{Sec7.3}
The practical value of the center-outward rank tests developed in the previous sections is illustrated with the following archeological application where classical methods fail to detect any treatment effect.  The data consist of $n=126$ measurements of MgO (Magnesium oxide), P$_2$O$_5$ (Phosphorus pentoxide), CoO (Cobalt monoxide), and~Sb$_2$O$_3$ (Antimony trioxide) (dimension $d=4$, thus) in natron glass vessels excavated  from three Syro-Palestinian sites in present-day Israel: Apollonia ($n_1=54$ observations),   Bet Eli'ezer ($n_2=17$ observations), and Egypt ($n_3=55$ observations); a fourth site only has two observations and   was dropped from the analysis. This dataset  has been   originally analyzed by~\cite{data-source} with the  objective  of detecting  possible differences among the three sites. 
  Bivariate plots of these four variables are shown in Figure~\ref{fig-data1}, where one can observe that the marginal distributions of~CoO, and Sb$_2$O$_3$ exhibit heavy tails and are very far from normal, and their  joint distribution     far from   elliptically symmetric. A traditional (pseudo-Gaussian) test here is  Pillai's trace test\footnote{Alternatives\vspace{-1mm} are Wilks' Lambda, the Lawley-Hotelling Trace, and Roy's largest root tests. In the two-sample case, they all coincide; else, they are   asymptotically equivalent.} reducing, in the two-sample case, to  Hotelling's classical $T$-square test.

   \begin{figure}[t!]
 \centering
 \includegraphics[width=0.8\textwidth, height=0.81\textwidth]{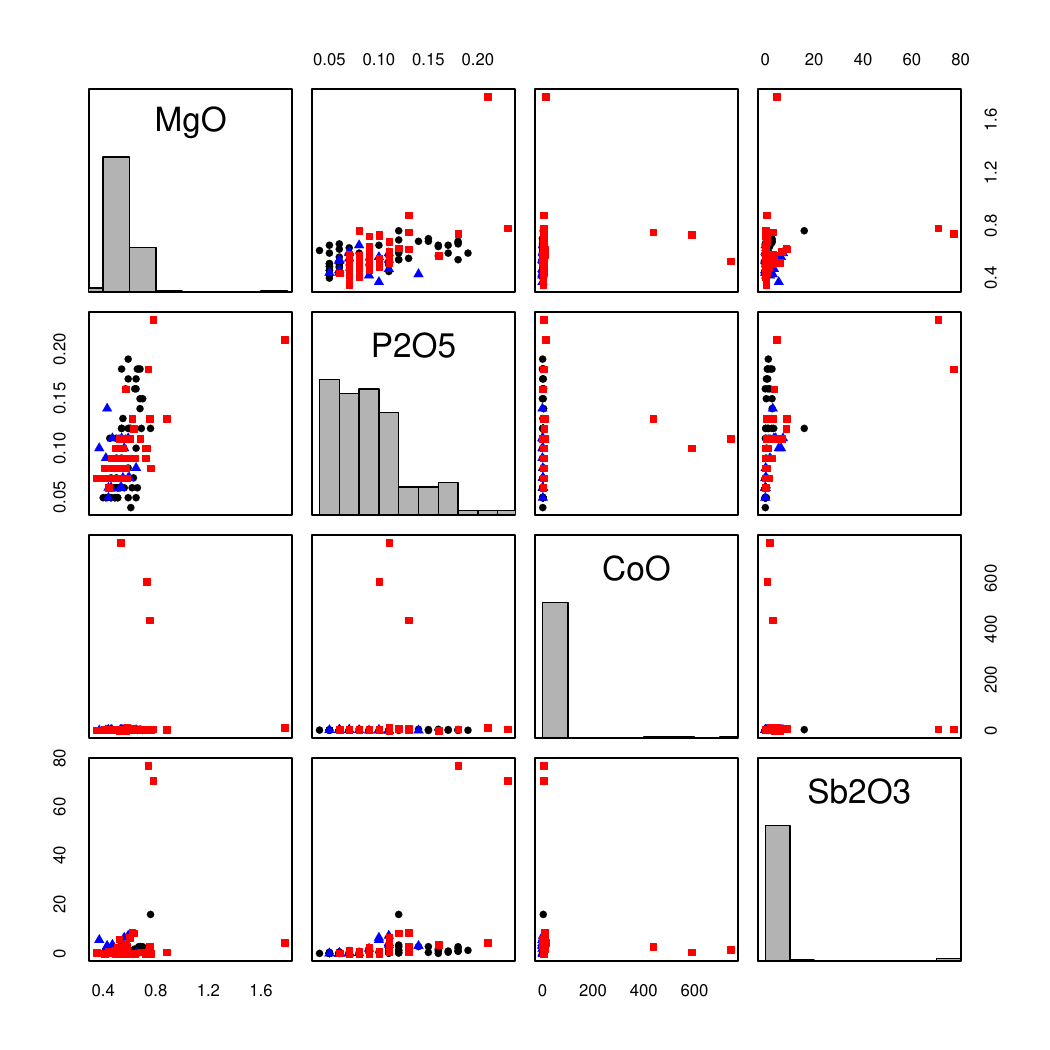}\vspace{-10mm}
 \caption{\small The content of  MgO, P$_2$O$_5$, CoO, and Sb$_2$O$_3$\vspace{-.0mm}  in natron glass vessels from Appo\-lonia (circles), Bet Eli'ezer (triangles), and Egypt (squares).
 }\label{fig-data1}
 \end{figure}

First, all the two-dimensional data subsets corresponding to the bivariate plots in Figure~\ref{fig-data1} were analyzed (six bivariate MANOVA models, thus).   Pillai's  test  yields  non-significant $p$-values for all  combinations, see Table~\ref{tab1}. But  the center-outward tests we are proposing in this paper do detect  significant differences between the three groups whenever the variable CoO is included in the analysis.  
Two versions of the center-outward ranks and signs are considered  in Table~\ref{tab1} below (c-o tests~I and  II, respectively; these two versions  correspond to two choices of the  grid~$\mathfrak{G}_n$,   with either $n_S=7$ and $n_R=18$  or~$n_S=18$ and~$n_R=7$---see Section~\ref{RSsec} for an explanation).  
\begin{table}[b!]
\centering
{
\begin{tabular}{llrrrr}
  \hline
&   & Pillai's test & Roy's test &c-o test I & c-o test II\\ 
  \hline
MgO & P$_2$O$_5$ &0.3547 & 0.3568 & 0.3817 & 0.0946 \\ 
  MgO & CoO & 0.1217 & 0.1592 & 0.0000 & 0.0000 \\ 
 MgO & Sb$_2$O$_3$ & 0.2268 & 0.3744 & 0.1865 & 0.3239 \\ 
P$_2$O$_5$ & CoO  & 0.1491 & 0.2747 & 0.0000 & 0.0000 \\ 
 P$_2$O$_5$ & Sb$_2$O$_3$  & 0.1957 & 0.3379 & 0.0569 & 0.2770 \\ 
 CoO & Sb$_2$O$_3$& 0.1453 & 0.1110 & 0.0000 & 0.0000 \\ 
   \hline
\end{tabular}}
\caption{\small $p$-values for the bivariate MANOVA  Pillai trace and   Wilcoxon center-outward  rank  tests\vspace{-.5mm}   based on~$n_R=7$, $n_S=18$  (c-o test I) and   $n_R=18$, $n_S=7$ (c-o test II), respectively.\vspace{-0mm}}
\label{tab1}
\end{table}

 Inspection of Table~\ref{tab1} reveals that, unlike Pillai's trace,  the Wilcoxon center-outward rank tests (c-o I and  II) reject the null hypothesis at significance level $\alpha=0.05$. 
As for the Wilcoxon tests based on elliptical ranks
(based on the sample covariance function), they  yield highly non-significant $p$-values for all couples of variables;  the corresponding results are not presented here. 
Next, the MANOVA comparison is conducted for the full $4$-dimensional dataset.  Pillai's and Roy's $p$-values are  $0.1553$ and  $0.2765$, respectively:  no difference  detected among the three groups, thus, at level \textcolor{blue}{$\alpha=0.05$}. In sharp contrast, the Wilcoxon center-outward rank test (with $n_R=7$ and $n_S=18$)  yields a $p$-value $10^{-15}$, which is highly significant.  The elliptical Wilcoxon rank test (based on the sample covariance matrix), on the other hand,  with  $p$-value $0.5827$, also  fails to detect anything at any  level~$\alpha\leq 0.5$.   

This, according to   archeological sources, might lead to revising some of  the conclusions made by \cite{data-source} on Middle-East economic exchanges between Egypt and Syro-Palestine in the {B}yzantine-{I}slamic transition period.

\section{Conclusion and perspectives} Classical multivariate analysis methods, which are daily practice in a number of applied domains, remain deeply marked by Gaussian and elliptical assumptions. In particular,  no distribution-free approach is available so far for hypothesis testing in multiple-output regression models, which include the fundamental two-sample and MANOVA models---except for the elliptical or Mahalanobis rank tests developed in \cite{HallinPain2005} which, however, require the strong assumption of elliptic symmetry---an assumption which is unlikely to hold in most applications. Based on the recent concept of center-outward ranks and signs, this paper proposes the first efficient fully distribution-free tests of the hypothesis of no treatment effect in that multiple-output  context, thereby extending to the multivariate case the classical H\' ajek  approach to univariate rank-based inference \citep{HajekSidak}. Simulations and an empirical example demonstrate the excellent performance of the method. This  lays the theoretical bases (asymptotic representation and asymptotic normality results for linear center-outward rank statistics) and  theoretical guidelines ({\it H\' ajek projection} of LAN central sequences) for the  development of a complete toolbox of distribution-free methods for multivariate analysis.%


\end{document}